\documentclass[reqno]{amsart}

\usepackage{amsthm, amssymb, latexsym, tabls, hhline,lscape}

\input xypic
\xyoption{all}
\xyoption{arc}

\newtheorem{prop}{Proposition}[section]
\newtheorem{thm}[prop]{Theorem}
\newtheorem{cor}[prop]{Corollary}
\newtheorem{lem}[prop]{Lemma}
\newtheorem{rem}[prop]{\it Remark}

\newcommand{\set}[1]{{\left\{ #1 \right\}}}  
\def\st{{\,|\,}}  
\def\implies{{\,\Rightarrow\,}} 

\def\std{{\textnormal{std}}}
\def\col{{\mathsf{col}}}

\def\UA{{\mathcal A}}

\def\UT{{\mathcal T}}

\def\ZM{{\mathbb Z}}
\def\NM{{\mathbb N}}

\def\SG{{\mathfrak S}}

\def\k{{\Bbbk}}
\def\b{{\mathsf b}}

\def\simsylv{\sim_{\textnormal{sylv}^*}}
\def\dotsimsylv{\dot\sim_{\textnormal{sylv}^*}}

\def\ipeak{{\mathring{\mathsf{peak}}}}
\def\peak{\mathsf{peak}}
\def\Des{\mathsf{Des}}

\def\IP{{\mathring P}}

\def\IPi{{\mathring\Pi}}
\def\ip{\mathring p}

\def\Comp{{\mathsf{N}^\ast}}

\def\IPeak{{\mathring{\mathcal P}}}
\def\Peak{{\mathcal P}}

\newcommand\QSym{\mathord{\mathrm{QSym}}}

\begin{document}

\title[Coloured peak algebras and Hopf algebras]{Coloured peak algebras and Hopf algebras}

\author{Nantel Bergeron}

\address{Department of Mathematics and Statistics,
York University, Toronto, Ontario M3J 1P3, Canada}

\email{bergeron@mathstat.yorku.ca}

\author{Christophe Hohlweg}

\address{The Fields Institute,
222 College Street,
Toronto, Ontario,
Canada M5T 3J1}

\email{chohlweg@fields.utoronto.ca}



\thanks{Bergeron is partially supported by NSERC and CRC, Canada}
\thanks{Hohlweg is partially supported by CRC}\maketitle

\begin{abstract}
For $G$ a finite abelian group, we study the properties of general
equivalence relations on $G_n=G^n\rtimes \SG_n$, the wreath
product of $G$ with the symmetric group $\SG_n$, also known as the
$G$-coloured symmetric group. We show that under certain
conditions, some equivalence relations give rise to subalgebras of
$\k G_n$ as well as graded connected Hopf subalgebras of
$\bigoplus_{n\ge o} \k G_n$. In particular we construct a
$G$-coloured peak subalgebra  of the Mantaci-Reutenauer algebra
(or $G$-coloured descent algebra). We show that the direct sum of
the $G$-coloured peak algebras is a Hopf algebra. We also have
similar results for a $G$-colouring of the Loday-Ronco Hopf
algebras of planar binary trees. For many of the equivalence
relations under study, we obtain a functor from the category of
finite abelian groups to the category of graded connected Hopf
algebras. We end our investigation by describing a Hopf
endomorphism of the $G$-coloured descent Hopf algebra whose image
is the $G$-coloured peak Hopf algebra. We outline a theory of
combinatorial $G$-coloured Hopf algebra for which the $G$-coloured
quasi-symmetric Hopf algebra and the graded dual to the
$G$-coloured peak Hopf algebra are central objects.
\end{abstract}

 \maketitle

\medskip
%
%

\section*{Introduction}

In recent work we have seen a vivid interest in uncovering
algebraic structures behind discrete objects and understanding the
relationship among them. These objects are particularly
interesting as their structure constants may encode various
invariants in geometry, in physics or in computer science.  This
is particularly true in \cite{aguiar bergeron nyman, bmws, hivert
novelli thibon 1, hivert novelli thibon 2, loday ronco, nyman,
reading, thibon}, just to cite a few, and in \cite{aguiar bergeron
sottile} we find the beginning of a general theory of
combinatorial Hopf algebras. Most of the algebraic structures
under study are subalgebras, subcoalgebras, Hopf subalgebras or
quotients of $\k[\SG]=\bigoplus_{n\ge 0} \k\SG_n$, where $\k$ is a
field (of characteristic $0$) and $\k\SG_n$ is the group algebra
of the symmetric group $\SG_n$. Besides the algebra structure
obtained from $\k\SG_n$, there are other algebraic structures on
$\k[\SG]$. For instance, Malvenuto-Reutenauer \cite{malvenuto
reutenauer} give to $\k[\SG]$ a structure of graded connected Hopf algebra.

We are interested in a natural generalization of $\k[\SG]$. Let
$G$ be a finite abelian group and let $G_n=G^n\rtimes \SG_n$ be
the wreath product of $G$ with $\SG_n$; this is viewed as a
$G$-colouring of $\SG_n$. Consider $\k[{\mathsf
G}]=\bigoplus_{n\ge 0}\k G_n$. In particular, when $G=\set{1}$ we
have that $\k[{\mathsf G}]=\k[\SG]$. In \cite{baumann hohlweg},
the authors extend to $\k[\mathsf G]$ the Hopf algebra structure
on $\k[\SG]$, 
this construction is functorial from the category of finite
abelian groups to the category of graded connected Hopf algebras.
It is then natural to ask if the multitude of substructures and
quotients of $\k[\SG]$ studied previously can be $G$-coloured. A
few answer have been given, see for instance \cite{aguiar bergeron
nyman, baumann hohlweg, novelli thibon, hivert novelli thibon3}.
In general this requires some work, and is not possible in all
cases. Here, our approach is to study properties of equivalence
relations on ${\mathsf G}=\bigoplus_{n\ge 0} G_n$ and to determine
when they give rise to algebraic structures. Many of the algebraic
objects obtained in this way are $G$-colouring of the structures
in the literature, such as the peak algebras \cite{aguiar bergeron
nyman, bmws, nyman} or the Loday-Ronco Hopf algebra of trees
\cite{loday ronco}. Our work is the unifying generalization of a
series of results starting in \cite{atkinson} and continuing in
\cite{baumann hohlweg, hivert novelli thibon, schocker}.

In Section~\ref{se:SymGroup}, we recall some known results about
Solomon's descent algebras and Nyman's peak algebras in symmetric
groups. We adopt the perspective of Atkinson \cite{atkinson} and
Schocker \cite{schocker} which is useful in subsequent sections.
In Section~\ref{se:ColouredPeaks} we recall the definition of the
Mantaci-Reutenauer algebra \cite{mantaci reutenauer}. These can be
viewed as $G$-coloured descent algebras \cite[\S5.1]{baumann
hohlweg}. For $G$, a finite abelian group, we develop a general
study of various equivalence relations on the $G$-coloured
symmetric groups. This allows us to determine which equivalence
relations give rise to subalgebras and when we have a Hopf algebra
structure (see Section~\ref{se:HopfAlgebras}). With this in hand,
we  show our first main theorem in \S\ref{ss:ColouredPeakAlg}:
that the span of coloured peak elements $\IPeak_n(G)$ form a
subalgebra of the Mantaci-Reutenauer algebra $\Sigma_n(G)$. We
then have a few important theorems in
Section~\ref{se:HopfAlgebras}  relating equivalence relations on
the $G$-coloured symmetric groups and certain graded connected
Hopf algebras. We also determine that under certain conditions,
this map is functorial from the category of finite abelian groups
to the category of graded connected Hopf algebras. This
generalizes related results in \cite{baumann hohlweg}. We give
some applications of this theory in Section~\ref{se:Appli}. In
particular, we describe a $G$-coloured peak Hopf algebra (a
$G$-colouring of the peak Hopf algebra from \cite{bmws}), and a
$G$-coloured Loday-Ronco Hopf algebra (a $G$-colouring of the
Loday-Ronco Hopf algebra \cite{loday ronco}). To our knowledge,
these two $G$-coloured algebras are new. We end this section by
giving examples of application of this theory in the
hyperoctahedral group ($\ZM/2\ZM$-permutations) with equivalence
relations not induced from symmetric groups. Finally, in
Section~\ref{se:Duality}, we introduce the $\Theta_G$ map. This is
a Hopf endomorphism of the $G$-coloured descent Hopf algebra whose
image is the $G$-coloured peak Hopf algebra. We then outline a
theory of combinatorial $G$-coloured Hopf algebras (generalizing
\cite{aguiar bergeron sottile}). We show that the Hopf algebra of
$G$-coloured quasisymmetric functions (the graded dual Hopf
algebra of the Mantaci-Reutenauer algebra, see \cite{baumann
hohlweg}) is the terminal object in the category of combinatorial
$G$-coloured Hopf algebras. We also state that the graded dual to
the $G$-coloured peak Hopf algebra is the terminal object in the
category of {\sl odd} combinatorial $G$-coloured Hopf algebras.

We wish to thank Hugh Thomas for fruitful discussions and Karin
Prochazka for her careful reading. The redaction of this article
was completed in May 2005, when the second author was visiting the
Institut Mittag-Leffler in Djursholm, Sweden.


\section{The descent algebra and the peak algebra in the symmetric group}
\label{se:SymGroup}

For $k,l\in \ZM$, we denote by $[k,l]$ the set of integers
$\set{k,k+1,\dots,l-1,l}$ if $k\leq l$ or $\set{l,l+1,\dots,k-1,k}$ if $l<k$.

\subsection{Words, Permutations and Symmetric groups}
\label{ss:PermutationsSymGroup}

Let us first recall some general notions about words and
permutations. An {\sl alphabet} $\UA$ is a totally ordered set.
Here we give the explicit description of the underlying order of
$\UA$ only when the context requires it. The elements of $\UA$ are
called {\sl letters}. We denote by $\UA^\ast$ the set of all words
$w=a_1 a_2\cdots a_n$ with letter $a_i\in\UA$ and $n\ge 0$. For
$n=0$ there is a unique empty word denoted by $(\,)$. The {\sl
length of a word} $w=a_1 a_2\cdots a_n$ is the integer $n$. An
{\sl injective word} is a word without repetition of letters. The
{\sl concatenation} $w\cdot w'$ of two words $w=a_1\dots a_n$ and
$w'=b_1\dots b_m$ in $\UA^\ast$ is the word $a_1\dots a_n b_1\dots
b_m$ of length $n+m$.

There are two important ways to envision the symmetric group
$\SG_n$. First as words encoding the bijections of the set $[1,n]$, and second as a Coxeter group.
We  use both points of view as needed. More specifically, let $\mathsf N=\set{1,2,3,\ldots
}$
be the alphabet with the usual order on integers.
A permutation $\sigma\in\SG_n$ is given as the injective  word
  $\sigma=\sigma(1)\sigma(2)\dots\sigma(n)$ in $\mathsf N^\ast$
  with letters in $[1,n]$.
Viewed as a Coxeter group, a permutation $\sigma\in\SG_n$ is now a word in $S_n^*$
where $S_n=\set{s_1,\dots,s_{n-1}}$ the
 set of simple transpositions $s_i = (i,i+1)$. It is well-known that $(\SG_n, S_n)$ is a Coxeter
 system of type $\mathsf A_{n-1}$ (e.g. \cite{humphreys}).
 We always denote by $1_n$ the identity of $\SG_n$.

 The {\sl length}
  of a permutation $\sigma\in\SG_n$ is
$$
\ell(\sigma)=\big|\set{(i,j) \st i<j\in [1,n]\textrm{ and }
\sigma(i)>\sigma(j )}\big|,
$$
namely, its {\sl number of inversions} (e.g.
\cite{humphreys}). We denote by $\omega_{n}$ the {\sl longest
element} in $\SG_n$.
This element is an involution and $\omega_n(i) =
n-i+1$, for all $i\in[1,n]$.

The {\sl standardization} of a word $w=a_1 a_2\cdots a_n$ of length $n$ in $\UA^\ast$, denoted
 by $\std(w)$,
 is the unique permutation $\sigma\in\SG_n$ such that for all $i<j$
 we have $\sigma(i)>\sigma(j)$ if and only if $a_i>a_j$.
For instance, for
$w=845912$ in $\mathsf N^\ast$ we have $\std(w)=534612$.

A {\sl composition} of $n$ is a word $\mathbf c$ in $\mathsf
N^\ast$ whose letters sum to $n$.
A composition $\mathbf c$ of $n$ is denoted by $\mathbf
c\vDash n$ and we also write $n=\|\mathbf c\|$.
To emphasize the fact that we have a composition, we write
$\mathbf c$ as a sequence $(c_1,\dots,c_k)$ of letter in $\mathsf
N$ instead of a concatenation of letters. There is a well-known
 bijection between compositions of $n$ and subsets of $[1,n-1]$
 defined by
\begin{eqnarray}
\label{eq:BijCompSubset1}%
 \mathbf c=(c_1,\dots, c_k)\vDash n &\mapsto&
I_{\mathbf c}=\{c_1,c_1+c_2,\ldots,c_1+c_2+\cdots+c_{k-1}\}.
\end{eqnarray}
The inverse of this map is
\begin{eqnarray}
\label{eq:BijCompSubset2}%
I=\set{i_1,i_2,\dots,i_k} &\mapsto& \mathbf c_I
=(i_1,i_2-i_1,\dots,i_k-i_{k-1},n-i_{k})\vDash n.
\end{eqnarray}

Let $\mathbf c=(c_1,c_2,\ldots,c_k)\vDash n$, we denote by
$\SG_{\mathbf c}$ the subgroup
 of $\SG_n$ generated by $S_{\mathbf c}=\set{s_i\st i\in [1,n]\setminus I_{\mathbf c}}$.
Such a subgroup is
 called a {\sl Young subgroup} (or  {\sl standard parabolic subgroup}). Set
$t_i=c_1+c_2+\cdots+c_i$ for all $i$. Given a $k$-tuple
$(\sigma_1, \sigma_2,\ldots,\sigma_k)\in\SG_{c_1}\times\SG_{c_2}
\times\cdots\times\SG_{c_k}$ of permutations, we define
$\sigma_1\times\sigma_2\times\cdots\times\sigma_k\in\SG_n$ as the
permutation that maps an element $a$ belonging to the interval
$[t_{i-1}+1,t_i]$ onto $t_{i-1}+\sigma_i(a-t_{i-1})$. This
assignment defines an isomorphism $\mathfrak
S_{c_1}\times\mathfrak S_{c_2}\times\cdots \times\mathfrak
S_{c_k}\simeq\SG_{\mathbf c}$. Finally, $\omega_{\mathbf c}
=\omega_{c_1}\times\dots\times\omega_{c_k}$ is the {\sl longest
element} in $\SG_{\mathbf c}$. For instance
$\omega_{(2,1,3)}=213654$.

For $\mathbf c=(c_1,\dots,c_k)$ and $t_i$ as above, the subset
$$X_{\mathbf c}=\bigl\{\sigma\in\mathfrak S_n\bigm|\forall i,\ \sigma
\text{ is increasing on the interval }[t_{i-1}+1,t_i]\bigr\}$$
is a system of representatives of the left cosets of $\mathfrak
S_{\mathbf c}$ in $\SG_n$. For example:
$$X_{(2,2)}=\{1234,1324,1423,2314,2413,3412\},\quad
\text{and}\quad X_{(\underbrace{1,1,\ldots,1}_{n\text{ times}})}
=\mathfrak S_n.$$
For $\sigma\in\SG_n$, there is a unique
 pair $(u,v)\in X_{\mathbf c}\times \SG_{\mathbf c}$, called
{the $\mathbf c$-components of $\sigma$}, such that
 $\sigma = uv$. Moreover $\ell(\sigma)=\ell(u)+\ell(v)$
(see \cite[\S1.10]{humphreys}). Specifically, write $\sigma$ as
 the concatenation  $\sigma_1\cdot\ldots\cdot\sigma_k$ of the words in $\mathsf N^\ast$
 such that  the length of the word $\sigma_i$ is $c_i$. It is then easy to check that
\begin{equation}
\label{eq:ParabolicCompenent}%
 v=\std(\sigma_1)\times
\std(\sigma_2)\times\cdots\times\std(\sigma_k) \in \SG_{\mathbf c}
\quad\text{and}\quad
u=\sigma v^{-1}\in
X_{\mathbf c}.
\end{equation}

Finally, let us recall   Deodhar's Lemma for $\SG_n$  \cite{deodhar}.

\begin{lem}
\label{le:Deodhar}
For  $\mathbf c\vDash n$, $x\in X_{\mathbf c}$ and  $s_i\in S_n$, either
\begin{enumerate}

\item[(i)] $s_i x\in X_{\mathbf c}$, or

\item[(ii)] $s_i x = xs_j$, with $s_j\in S_{\mathbf c}$.
In this case,  $x(j)=i$ and {{$x(j+1)=i+1$}}.

\end{enumerate}
\end{lem}

\subsection{Descent sets and Solomon's descent algebra}
\label{ss:DescentSets}

For $\sigma\in\SG_n$, the {\sl descent set} of $\sigma$ is
$$
\Des(\sigma)=\set{i\in[1,n-1]\st \sigma(i)>\sigma(i+1)}.
$$
The {\sl descent composition} $D(\sigma)$ of $\sigma\in\SG_n$ is
defined via Eq.~(\ref{eq:BijCompSubset2}): $D(\sigma)=\mathbf c_{\Des
(\sigma)}$. For any $I\subseteq [1,n-1]$ we
have $\Des(\omega_{\mathbf c_{[1,n-1]\setminus I}})=I$,
in particular all compositions of $n$ are  descent compositions.

Denote by $\ZM\SG_n$ the group algebra of $\SG_n$. For each
 $\mathbf c\vDash n$, we set
$
d_{\mathbf c} =\sum_{D(\sigma)=\mathbf c} \sigma.
$
Then $$ \Sigma_n=\bigoplus_{\mathbf c\vDash n} \ZM d_{\mathbf c}
$$
is a subalgebra of $\ZM\SG_n$, called the Solomon descent algebra
\cite{solomon}. Let us recall the ingredients of Atkinson's \cite{atkinson} alternate proof of
this result (see also \cite[\S4]{schocker}).  Let $s_i\in S_n$ and $\sigma\in\SG_n$.
We  define an equivalence relation as follows:
\begin{equation}
\label{eq:EquivDescentSn} \sigma\, \dot\sim_D\, s_i\sigma\quad \iff\quad
\sigma^{-1}s_i\sigma\notin S_n\quad\iff\quad
|\sigma^{-1}(i)-\sigma^{-1}(i+1)|>1 .
\end{equation}
It is easily seen that $\dot\sim_D$ is a symmetric relation. The
reflexive and transitive closure of $\dot\sim_D$ is called
\textit{the descent equivalence relation} and is denoted by
$\sim_D$.

\begin{prop}
\label{pr:AtkinsonDescent} $\quad$
\begin{enumerate}
\item[(i)] For $s_i\in S_n$ and $\sigma\in\SG_n$, if
$s_i\sigma\dot\sim_D \sigma$, then $D(s_i \sigma)=D(\sigma)$.

\item[(ii)] For $\sigma,\tau\in\SG_n$,\quad $D(\sigma)=D(\tau)$ if
and only if $\sigma\sim_D \tau$.
\end{enumerate}
\end{prop}

\begin{rem}
\label{re:DescentLeftConnected} \textnormal{This proposition was
first stated by Tits \cite[Theorem 2.19]{tits} in a more general
context (see also \cite{atkinson, bjorner wachs}).}

\textnormal{ Recall that a subset $X$ of $\SG_n$ is {\sl
left-connected} if for each $\sigma,\tau\in X$, there is a
sequence $\sigma_1 = \sigma , \sigma_2,\dots ,\sigma_k = \tau$ of
elements in $X$ such that $\sigma_{j+1}\sigma_{j}^{-1}\in S_n$,
for all $j\in[1,k-1]$. That is, these classes can be seen as a set
of adjacent nodes in the {\sl permutohedron}. In particular,
Proposition~\ref{pr:AtkinsonDescent} shows that the sets of
permutations whose descent composition is fixed are
left-connected.}
\end{rem}

\begin{cor}
\label{co:DescentSets} Let $I\subseteq [1,n-1]$. If there is $j\in
[1,n-1]$ such that $j-1,j\in I$ and $j+1\notin I$, then there is
$\sigma\in\SG_n$ such that $\Des(\sigma)=I$, $\sigma(j)=2$ and
$\sigma(j+1)=1$.
\end{cor}
\begin{proof}
Fix $I$ and $j$ as above and consider the set
  $$T=\set{\tau\in\SG_n \st \Des(\tau)=I,\, \tau(j)=\tau(j+1)+1,\, \hbox{ and }\forall i<\tau(j+1),\,
  \tau^{-1}(i)<j
    }.$$
Let $\mathbf c$ be such that $\Des(\omega_{\mathbf c})=I$. Clearly
$\omega_{\mathbf c}\in T$, so the set $T$ is nonempty. Start with
any $\tau\in T$ and set $k=\tau(j+1)$.  If $k=1$, then we have found the desired
$\sigma=\tau$. If $k>1$, then let $\tau_1=s_{k-1}\tau$. Since
$\tau\in T$,  $\tau^{-1}(k-1)<j<\tau^{-1}(k)$. Therefore
Eq.~(\ref{eq:EquivDescentSn}) implies that $\tau_1\sim_D \tau$.
From Proposition~\ref{pr:AtkinsonDescent} we have that
$\Des(\tau_1)=I$. The fact that $j-1\in I$ gives that
$\tau_1^{-1}(k)=\tau^{-1}(k-1)<j-1$, so
$\tau_1^{-1}(k+1)=\tau^{-1}(k+1)=j>j-1>\tau_1^{-1}(k)$ and again
we can apply the equivalence relation to get
$\tau_2=s_k\tau_1\sim_D\tau_1\sim_D\tau$. We have that $\tau_2\in
T$ with $\tau_2(j+1)=k-1$. We can repeat this process $k-1$ times
to get the desired $\sigma$ with $k=1$.
\end{proof}

\noindent For $\sigma\in\SG_n$ and $\mathbf c,\mathbf d$
compositions of $n$, we define
$$
D(\mathbf c,\mathbf d,\sigma) = \set{(u,v) \in \SG_n \times \SG_n
\st D(u)=\mathbf c,\ D(v)=\mathbf d,\ uv=\sigma }.
$$
Let us fix $s_i \in S_n$ and  let $(u,v) \in \SG_n \times \SG_n$.
If $u\, \dot\sim_D\, s_i u$, we set $\psi^D_i(u,v)=(s_i u,v)$. If
$u \not\sim_D s_i u$, then we set $\psi^D_i (u,v)=(u,u^{-1}s_i u
v)$. Note that in the second case we have $u^{-1}s_i u \in S_n$. Since
$(\psi^D_i)^2 = \text{Id}_{\SG_n \times \SG_n}$, we have that
$\psi^D_i$ is a bijection.

\begin{lem}[Atkinson~\cite{atkinson}]\label{le:Psi}
If $\sigma \in \SG_n$ and $s_i \in S_n$ are such that $\sigma
\dot\sim_D s_i \sigma$, then $\psi^D_i\big(D(\mathbf c,\mathbf
d,\sigma)\big)=D(\mathbf c,\mathbf d,s_i\sigma)$, for all $\mathbf
c,\mathbf d$ compositions of $n$.
\end{lem}

Lemma~\ref{le:Psi} and
Proposition~\ref{pr:AtkinsonDescent} show that $\Sigma_n$ is a
subalgebra of $\ZM\SG_n$. Indeed,
  $|D(\mathbf c,\mathbf d,\sigma)|$ depends only on  $\mathbf e=D(\sigma)$. Setting
 $a_{\mathbf c,\mathbf d,\mathbf e}=|D(\mathbf c,\mathbf d,\sigma)|$, we have
$$
d_{\mathbf c} d_{\mathbf d} =\sum_{\mathbf e\vDash n} a_{\mathbf
c,\mathbf d,\mathbf e} d_{\mathbf e} .
$$

\subsection{Nyman's peak algebra}
\label{ss:IntPeaks}

Let $\sigma\in\SG_n$, the {\sl peak set} of $\sigma$ is
$$
\ipeak(\sigma)=\set{i\in[1,n-1]\st
\sigma(i-1)<\sigma(i)>\sigma(i+1),\ i>1}.
$$
This set is sometimes called the {\sl interior peak set}
\cite{aguiar bergeron nyman,aguiar nyman orellana}. Contrary
to descent sets, all the subsets of $[1,n-1]$ are not peak sets.
In fact $I\subseteq [1,n-1]$ is a peak set if and only if $I$ satisfies
the condition: if $i\in I$, then $i\ge 2$ and $i-1\notin I$. The {\sl  peak
composition} $\IP(\sigma)$ of $\sigma\in\SG_n$ is defined by
Eq.~(\ref{eq:BijCompSubset2}): $\IP(\sigma)=\mathbf
c_{\ipeak(\sigma)}$. Denote by $\IPi_n$
 the set of all compositions of $n$ which are peak compositions.
  {}From the above discussion on peak sets, it is  obvious that
\begin{equation}
\label{eq:DescripInteriorPeakComp}%
\IPi_n=\set{\mathbf (c_1,\dots,c_k)\vDash n \st
c_i>1 \hbox{ for } 1\le i\le k-1}.
\end{equation}
Moreover, it is clear
 that for all $\sigma,\tau\in\SG_n$
\begin{equation}
\label{eq:InteriorPeaksInDescent}
\ipeak(\sigma)\subseteq\Des(\sigma) \quad\textrm{and}\quad
D(\sigma)=D(\tau)\ \implies\ \IP(\sigma)=\IP(\tau).
\end{equation}

 For each $\mathbf c\in\IPi_n$, we set
$
\ip_{\mathbf c} =\sum_{\IP(\sigma)=\mathbf c} \sigma.
$
  Nyman \cite{nyman} has shown that
$$
\IPeak_n=\bigoplus_{\mathbf c\in\IPi_n} \ZM \ip_{\mathbf c}
$$
is a subalgebra of $\Sigma_n$. We call $\IPeak_n$ the {\sl peak
algebra}. Adapting the argument of Atkinson, Schocker
\cite{schocker} gives a new proof of Nyman's result. We recall here the main  ingrediants.
 Let $s_i\in S_n$ and $\sigma\in\SG_n$. Then
\begin{eqnarray}
\label{eq:EquivIntPeakSn}
\sigma \dot\sim_\IP s_i\sigma &\iff&
s_i\sigma\dot\sim_D \sigma
\textrm{ or } i=1
\\
&\iff&
|\sigma^{-1}(i)-\sigma^{-1}(i+1)|>1
\textrm{ or }i=1\,.
\end{eqnarray}
 The reflexive and transitive closure of $\dot\sim_\IP$ is called
\textit{the  peak equivalence}, and is denoted by $\sim_\IP$.
\begin{prop}[{Schocker~\cite{schocker}}]
\label{pr:IntPeakSets} $\quad$
\begin{enumerate}
\item[(i)] For $\sigma\in\SG_n$ and $s_i\in S_n$, if
$s_i\sigma\dot\sim_\IP \sigma$, then $\IP(s_i\sigma)=\IP(\sigma)$.

\item[(ii)] For $\sigma,\tau\in\SG_n$, \quad
$\IP(\sigma)=\IP(\tau)$ if and only if $\sigma \sim_\IP \tau$.
\end{enumerate}
\end{prop}

\begin{rem}\label{re:IntPeakConnected}
\textnormal{As a consequence, the sets of permutations  having the
same peak compositions are left-connected.}

\textnormal{The proof of the above proposition uses
Corollary~\ref{co:DescentSets}. We  follow a similar recipe
for the  proof of Proposition~\ref{pr:PeakSetsEquiva} in
\S\ref{ss:PeakSets}.}
\end{rem}

 For $\sigma\in\SG_n$ and $\mathbf c,\mathbf d\in\IPi_n$, we set
$$
\IP(\mathbf c,\mathbf d,\sigma) = \set{(u,v) \in \SG_n \times
\SG_n \st \IP(u)=\mathbf c,\ \IP(v)=\mathbf d,\ uv=\sigma }.
$$
Fix $s_i \in S_n$ and  let $(u,v) \in \SG_n \times
\SG_n$. Then
$$
\psi_i^\IP (u,v) = \left\{\begin{array}{ll}
(s_1 u,v)&\textrm{if }i=1\\ \\
\psi^D_i (u,v)&\textrm{otherwise}
\end{array}\right.
$$
is an involution on $\SG_n \times \SG_n$. The following Lemma,
Eq.~(\ref{eq:InteriorPeaksInDescent}), and
Proposition~\ref{pr:IntPeakSets} show that $\IPeak_n$ is a
(non-unitary) subalgebra of $\Sigma_n$.

\begin{lem}[{Schocker~\cite{schocker}}]
\label{le:PsiIntPeak} If $\sigma \in \SG_n$ and $s_i \in S_n$
are such that $\sigma \dot\sim_\IP s_i \sigma$, then
$\psi^\IP_i\big(\IP(\mathbf c,\mathbf d,\sigma)\big)=\IP(\mathbf
c,\mathbf d,s_i\sigma)$, for all $\mathbf c,\mathbf d \in \IPi_n$.
\end{lem}

%
%

\section{Coloured peak algebras}
\label{se:ColouredPeaks}

\subsection{Coloured words}
\label{ss:Coloured words}

Let $G$ be a finite abelian group. We denote by $1^G$ the identity
of $G$. We also refer to $G$ as a set of {\sl colours}.

 Let $\UA$ be an alphabet. We denote  by $\UA^G$ the alphabet $\mathcal A\times G$,
and $a^g$ a letter $(a,g)$ in $\mathcal A^G$. We say that $\UA$ is
 {\sl coloured} by $G$.  The {\sl colour} of the letter $a^g$, denoted by
 $\col(a^g)$, is $g$.
Let $\mathbf g=(g_1,\dots,g_n)$ be a sequence of elements in $G$
 and $w=a_1 a_2 \dots a_n$ be a word in $\UA^\ast$.
 We denote by $w^\mathbf{g}$ the word
 $a_1^{g_1}\cdot a_2^{g_2}\cdot\ldots\cdot  a_n^{g_n}$ in $(\UA^G)^\ast$.
If $g_1=g_2=\dots=g_n=g$, we simply write
 $w^{\mathbf g}=w^g$. Therefore we can identify  $\big((\UA^\ast)^G\big)^\ast$ and $(\UA^G)^\ast$.
The {\sl absolute value map} $|\, .\,|\colon(\UA^G)^\ast\to
\UA^\ast$ maps a word $w^{\mathbf g}$ to the word $w$ in
$\UA^\ast$. This map satisfies $|w^{\mathbf g}\cdot u^{\mathbf
h}|=|w^{\mathbf g}|\cdot| u^{\mathbf h}|$.

A {\sl $G$-composition} $\mathbf c^{\mathbf g}$ is a word
in $(\mathsf N^G)^\ast$, and $|\mathbf c^{\mathbf g}|=\mathbf c$
 is a composition.
A $G$-composition of $n$, denoted by $\mathbf c^{\mathbf g} \vDash_G n$,
is a $G$-composition such that
$\|\mathbf c^{\mathbf g}\|=\|\mathbf c\|=n$.

 The {\sl rainbow decomposition} of a word $w^{\mathbf g}$ in $(\UA^\ast)^G$ is
 the unique sequence of non-empty words $w_1^{g_1}$, $w_2^{g_2}$,
\dots, $w_k^{g_k}$
 such that for all $i\in[1,k-1]$, $g_i\not = g_{i+1}$, and
\begin{equation}
\label{eq:RainbowDecomp}
w^{\mathbf g}= w_1^{g_1}\cdot w_2^{g_2}\cdot \ldots \cdot w_k^{g_k}
\end{equation}
Let $c_i$ be the length of the word $w_i$. The {\sl rainbow composition}
 of $w^{\mathbf g}$ is $\mathbf r(w)=(c_1^{g_1},c_2^{g_2},\dots,c_k^{g_k})$.
For instance, if $g\not= h\in G$ and
$w^{\mathbf g}=1^g 4^g 2^g 6^h 7^g 3^h 5^h$,
 then  the rainbow decomposition and the rainbow composition are
$$
 w^{\mathbf g}=142^g \cdot 6^h \cdot 7^g \cdot 35^h \quad\textrm{and}
\quad\mathbf r(w^\mathbf{g})=(3^g,1^h,1^g,2^h).
$$
  This procedure is the first part
 of the construction of the {\sl $G$-descent composition} $D_G (w^\mathbf{g})$
given by Mantaci and Reutenauer \cite{mantaci reutenauer}.
The {\sl descent composition
$D(w)$}
 of a word in $\UA^\ast$ is
\begin{equation*}
\label{eq:DescentWords}
D(w)=D\big( \std(w)\big).
\end{equation*}
If the length of $w$ is $n$, then $D(w)\vDash n$.
 The {\sl $G$-descent composition} $D_G(w^\mathbf{g})$ of a word
 $w^\mathbf{g}$  in $(\UA^G)^\ast$ is obtained as follows. Let
 $w^{\mathbf g}=w_1^{g_1}\cdot w_2^{g_2}\cdot \ldots \cdot w_k^{g_k}$
 be the rainbow decomposition of $w^\mathbf{g}$. We have
\begin{equation}
\label{eq:GDescentWords} D_G(w^{\mathbf g})=D(w_1)^{g_1}\cdot
D(w_2)^{g_2}\cdot \ldots\cdot D(w_k)^{g_k} .
\end{equation}
Again, if the length of $w^{\mathbf g}$ is $n$,
then $D_G(w^{\mathbf g})\vDash_G n$. For instance, with $g\not= h\in G$ and
$w^{\mathbf g}=1^g 4^g 2^g 6^h 7^g 3^h 5^h$, we have
$D_G(w^\mathbf{g})=(2^g,1^g,1^h,1^g,2^h)$.

\subsection{Coloured permutations and the Mantaci-Reutenauer algebra}
\label{ss:MRAlgebra}

We denote by $G_n=G\wr \SG_n=G^n\rtimes \SG_n$ the wreath product of
 $G$  with $\SG_n$. An element $w\in
G_n$ is the product $\alpha=\sigma.(g_1,g_2,\dots , g_n)$ where
$\sigma \in \SG_n$ and $g_i\in G$. We call {\sl coloured
permutations} the elements
 of $G_n$.   The multiplication in $G_n$ comes from
 the following commutation rule between elements of $\SG_n$ and elements of  $G^n$:
\begin{equation}
\label{eq:CommRule} \sigma.(g_1,g_2,\dots,g_n) =
(g_{\sigma^{-1}(1)},g_{\sigma^{-1}(2)},\dots,g_{\sigma^{-1}(n)}
).\sigma.
\end{equation}
 The subgroup of $G_n$ consisting of the element
 $\sigma.(1^G,\dots , 1^G)$ is isomorphic to $\SG_n$. In this
 case, we  write $\sigma$ instead of $\sigma.(1^G,\dots , 1^G)$. Therefore $1_n$ is the identity of $G_n$.
 For $\mathbf g=(g_1,g_2,\dots , g_n)$,
instead of the multiplicative notation $\sigma.\mathbf g$, it is more convenient to  use the {\it word notation}
$\sigma^{\mathbf g}$ in $(\mathsf N^\ast)^G$.

A subset $X$ of $G_n$ is {\sl left-connected} if for each $\alpha,\beta\in X$,
there is a sequence $\alpha=x_1 , x_2,\dots ,x_k
= \beta$ of elements in $X$ such that
$x_{j+1}x_{j}^{-1}\in S_n$, for all $j\in[1,k-1]$. In particular, all the
 elements in a left-connected set have the same sequence of colours.

Denote by $G_{n,m}$ the subgroup of $G_{n+m}$   isomorphic to $G_n
\times G_m$, obtained by the map
\begin{eqnarray*}
G_n\times G_{m} & \to & G_{n+m}\\
(\sigma^\mathbf{g},\tau^\mathbf{h})&\mapsto&
(\sigma\times\tau)^{\mathbf g\cdot\mathbf h}.
\end{eqnarray*}
For any $\alpha\in G_n$ we obtain a unique decomposition
$\alpha=\sigma^\mathbf{g}=u v^\mathbf{g}$, where
$(u,v^\mathbf{g})\in X_{n,m}\times G_{n,m}$ and $(u,v)$ are the
$(n,m)$-components of $\sigma$ (see
\S\ref{ss:PermutationsSymGroup}).

 Denote by $\ZM G_n$ the group algebra of $G_n$. For each $G$-composition
 $\mathbf c^{\mathbf g}$ of $n$, we set
$
d^G_{\mathbf c^{\mathbf g}} =\sum_{D_G(\alpha)=\mathbf c^{\mathbf
g}} \alpha.
$
Then
$$
\Sigma_n (G)=\bigoplus_{\mathbf c^{\mathbf g}\vDash_G n} \ZM
d^G_{\mathbf c^{\mathbf g}}
$$
is a subalgebra of $\ZM G_n$, called the {\sl Mantaci-Reutenauer
algebra} \cite{mantaci reutenauer}. It is clear,
$\Sigma_n(\set{1}) = \Sigma_n$. The Atkinson proof using an
equivalence relation on $\SG_n$ has recently been extended to
$G_n$ for $\Sigma_n(G)$ \cite{baumann hohlweg}.
In the subsections \ref{ss;GradedEquivRel}  to \ref{ss:MultipkicativeProperty},
we develop a general theory of equivalence relations
for the coloured symmetric groups.
This will allow us to introduce  {\sl coloured peak algebras} in $\ZM G_n$
and will be useful in  \S\ref{se:HopfAlgebras} to study Hopf structures.

\subsection{Graded connected equivalence relations}
\label{ss;GradedEquivRel}

Recall that a {\sl graded set} is a pair $(\mathsf E,\deg)$
 where $\mathsf E$ is a set and $\deg\colon\mathsf E \to \NM$ is a map. For simplicity,
  we  write $\mathsf E$ instead
 of $(\mathsf E,\deg)$ whenever possible. For $n\in
\NM$, we set $E_n=\deg^{-1}\set{n}$. Then
$$
\mathsf E=\bigoplus_{n \in \NM} E_n .
$$
For instance, $\mathsf G=\bigoplus_{n\in\NM} G_n$,
$\SG=\bigoplus_{n\in\NM}\SG_n$ and
$(\mathsf N^G)^\ast$ with the length of words are graded sets. On $(\mathsf N^G)^\ast$ we can also consider
$\big((\mathsf N^G)^\ast,\|\cdot\|\big)$  which gives a different graded set structure.

A {\sl graded map} $\rho$ is a map $\rho\colon\mathsf E' \to \mathsf E$
between two graded sets such that  $\rho(E'_n)$ is contained in $E_n$. A {\sl graded equivalence relation} on a graded set $\mathsf E$ is an equivalence relation on each $E_n$.
Given  $\rho\colon\mathsf E' \to \mathsf E$ a graded map, we
  define a graded equivalence relation $\sim_\rho$ on $\mathsf E'$
  where $e'_1\sim_\rho e'_2$ if and only if
 $\rho(e'_1)=\rho(e'_2)$. In this case, we can view $\rho(E'_n)$
as a parametrization of the set of equivalence
 classes of $E'_n$. Conversely, to any graded equivalence relation $\sim$
  we can associate a surjective graded map $\rho$ such that $\sim$ is $\sim_\rho$.


Let $\rho\colon\mathsf G \to \mathsf E$ be a graded map and $\mathbf
c=(c_1,\dots,c_k)\vDash n$. We use $\rho$ to define a  map on
$G_\mathbf{c}$ as follows:
\begin{eqnarray}
\label{eq:RhoParabolic}
 \rho_\mathbf{c}\colon G_\mathbf{c}&\to& \underbrace{\mathsf E\times\mathsf E\times \ldots\times\mathsf E}_{k\textrm{
times}}\\
v_1\times v_2\times\dots\times v_k&\mapsto&
\rho(v_1)\times\rho(v_2)\times\dots\times\rho(v_k)\,.
\end{eqnarray}
It is clear that $\rho_\mathbf{c}=\rho_{c_1}\times\rho_{c_2}\times\cdots
\times\rho_{c_k}$ and $\rho_{c_i}=\rho$ for all
$i\in[1,k]$. This multiplicative notation  is useful for
looking at induction from Young subgroups (see
\S\ref{ss:InductionProperty}).

A {\sl  graded connected map on $\mathsf G$} is a graded
 map $\rho\colon\mathsf G \to \mathsf E$ where each equivalence class under $\sim_\rho$
  is left-connected. In particular $\sim_\rho$ is the
 transitive and reflexive closure of the symmetric relation
 $\dot\sim_\rho$ defined as follows: for $n\in\NM$, $w\in G_n$
and $s_i\in S_n$ we have
$$
w\dot\sim_\rho s_i w \iff \rho(w)=\rho(s_i w).
$$
In this case, we say that $\sim_\rho$ is a {\sl graded connected
equivalence relation on $\mathsf G$}.
 For instance $D,\IP\colon\SG\to (\Comp,\|\,.\,\|)$ are graded connected maps on ${\mathsf G}=\SG$ (when $G=\set{1}$).

\subsection{The induced graded connected maps $\rho_G$ and $\rho^G$}
\label{ss:GenEquivG}

Let $\rho\colon\SG\to\mathsf E$ be a fixed graded connected map on $\SG$. For $G$  a finite abelian group we  define a graded connected
map on $\mathsf G$
as follow.
Let $n\in\NM$, $\sigma^{\mathbf g}\in G_n$ and $s_i\in S_n$.
First, identify in the rainbow decomposition $\sigma^{\mathbf
g}=\sigma_1^{g_1}\dots\sigma_k^{g_k}$ the subword $\sigma_{\nu}$
containing $i$ as a letter and the subword $\sigma_\mu$ containing
$i+1$ as a letter. We then define
\begin{equation}
\label{eq:CongEquivG}%
\sigma^{\mathbf g}\dot\sim_\rho^G s_i \sigma^{\mathbf g} \iff
\left\{
\begin{array}{l}
\nu\not=\mu\quad
\textrm{or}\\
 \nu=\mu\quad\textrm{and}\quad \rho(s_i\sigma)=\rho(\sigma)\,.
\end{array}\right.
\end{equation}
Denote by $\sim_\rho^G$ the transitive and reflexive  closure of
this relation. We let $E_n(G)=G_n\big/\! \sim_\rho^G$ be the set of equivalence classes
on $G_n$ and define $\rho_G$ the projection
\begin{equation}
\label{eq:CongMapG} %
\rho_G\colon\mathsf G\to \mathsf E(G)=\bigoplus_{n\geq 0} E_n(G)\,.
\end{equation}
This is a graded connected map on $\mathsf G$.

\begin{rem}\label{re:GenerGradedConnectedMap}
\textnormal{The relations $\dot\sim_D^{\set{1}}$ and
$\dot\sim_\IP^{\set{1}}$ are not the relations $\dot\sim_D$ and
$\dot\sim_\IP$ as defined in \S\ref{se:SymGroup}, but both induce
the same equivalence relations. For instance\break
$1423\dot\sim_\IP^{\set{1}}1432$  and $1423\sim_\IP1432$, but
$1423$ and $1432$ are not related by $\dot\sim_\IP$.  More
generally if $\dot\sim_\rho$ is a graded symmetric relation on
$\SG$ whose transitive and reflexive closure is $\sim_\rho$, we
can define a symmetric graded relation $\dot\sim_{\rho_G}$ on
$\mathsf G$ as in \S\ref{se:SymGroup}, namely:
\begin{equation}
\label{eq:CongEquivGBIS}%
\sigma^{\mathbf g}\dot\sim_{\rho_G} s_i \sigma^{\mathbf g} \iff
\left\{
\begin{array}{l}
\nu\not=\mu\quad
\textrm{or}\\
 \nu=\mu\quad\textrm{and}\quad s_i\sigma\dot\sim_\rho\sigma\,.
\end{array}\right.
\end{equation}
The transitive and reflexive closure of
$\dot\sim_{\rho_G}$ is precisely $\sim_\rho^G$. In fact
$$
\sigma^{\mathbf g}\dot\sim_{\rho_G} s_i \sigma^{\mathbf g}
\implies \sigma^{\mathbf g}\dot\sim_{\rho}^G s_i \sigma^{\mathbf
g},
$$
but the converse is not true. We will use Eq.~(\ref{eq:CongEquivG})
  for studying general properties of graded connected maps on $\mathsf G$ and Eq.~(\ref{eq:CongEquivGBIS})
 for studying particular cases of graded connected maps like
 the $G$-descent composition map or the $G$-peak composition map.}
\end{rem}

We have another way to define a graded map on $\mathsf G$ from $\rho\colon\SG\to\mathsf E$.
Again, let $\rho\colon\SG\to\mathsf E$ be a fixed graded connected map on $\SG$.
We first extend $\rho$ to a graded map $\mathsf N^\ast\to\mathsf E$ by setting $\rho(w) = \rho\big(\std(w)\big)$. Then, for
 $\sigma^{\mathbf g}=\sigma_1^{g_1}\dots\sigma_k^{g_k}$
 the rainbow decomposition of $\sigma^{\mathbf g}\in G_n$ we set
 \begin{equation}
\label{eq:DefRhoG} \rho^G(\sigma^{\mathbf g})=
 \rho (\sigma_1)^{g_1}\cdot\rho(\sigma_2)^{g_2}\cdot\ldots\cdot \rho (\sigma_k)^{g_k}.
\end{equation}
This defines a graded connected map $\rho^G\colon{\mathsf G}\to ({\mathsf E}^G)^\ast$.
Observe that if $\mathsf E$ is the set of words $\UA^\ast$ on an
alphabet $\UA$, then $
(\mathsf E^G)^\ast = ((\UA^\ast)^G)^\ast=(\UA^G)^\ast.$

The graded equivalence relation $\sim_{\rho^G}$ can be obtained as the transitive and reflective closure of the following graded connected symmetric relation $\dot\sim_{\rho^G}$ on $\mathsf G$.
For this, let $\sigma^{\mathbf g}\in G_n$ and $s_i\in S_n$, and select $\mu$ and $\nu$ as in Eq.~(\ref{eq:CongEquivG}).
\begin{equation}
\label{eq:GenEquivG}%
\sigma^\mathbf{g}\dot\sim_{\rho^G} s_i \sigma^\mathbf{g} \iff
\left\{
\begin{array}{l}
\nu\not=\mu\quad
\textrm{or}\\
 \nu=\mu\quad\textrm{and}\quad
\rho(\sigma_\nu)=\rho(s_i \sigma_\nu ).
\end{array}\right.
\end{equation}


\begin{rem}\label{re:SurjGradedMap} \textnormal{If $\rho$ is surjective, then $\rho^G$ is
surjective. In particular, from
Eq.~(\ref{eq:DefRhoG}) we get a bijection
\begin{eqnarray}
\label{eq:BijInducedEquiv} (E^G)^\ast_n&\to&
\bigoplus_{{\scriptstyle (c_1,\dots,c_k)\vDash n \atop
\scriptstyle
g_i\not = g_{i+1} \in G}} %
E_{c_1}^{g_1}\times\dots\times E_{c_k}^{g_k}.
\end{eqnarray}
Thus
\begin{equation}
\label{eq:CardinalInducedEquiv}
|(E^G)^\ast_n|=|G|\sum_{(c_1,\dots , c_k)\vDash n} \left(
|G|-1\right)^{k-1} |E_{c_1}|\dots |E_{c_k}|.
\end{equation}
} \end{rem}

In general, the graded maps $\rho^G$ and $\rho_G$ have nothing in
common. For
 instance, let $\sim_\rho$ be defined on $\SG_n$ by
$s_i\sigma \dot\sim_\rho \sigma$ if $\sigma^{-1}([i,i+1])=[1,2]$.
If $g\not =h$  are two elements in $G$,  then
$1^{g}2^{g}3^{h}4^{h}\dot\sim_{\rho^G} 1^{g}2^{g}4^{h}3^{h}$.
 But these two elements are in distinct equivalence classes
under $\rho_G$.
 The map $\rho^G$  is easier to study but does not in general  induce the
properties we wish to extend (subalgebra, Hopf subalgebra, etc).  However, if $\rho$ possesses  the {\sl induction property} (\S\ref{ss:InductionProperty})
and the {\sl freeness property} (\S\ref{ss:MultipkicativeProperty}), the two relations  $\sim_{\rho_G}$ and  $\sim_{\rho^G}$ will coincide.

\subsection{Induction Property}
\label{ss:InductionProperty}

We say that a graded map $\rho\colon\SG\to\mathsf E$ or  its graded
equivalence relation $\sim_\rho$ have {\sl the induction property}
if the following condition is satisfied:

\begin{itemize}
\item[(IP)] For any  $n,m\in\NM$ and $(e_n,e_m)\in E_n\times E_m$,
 there is $X\subseteq E_{n+m}$ such that
$$
X_{(n,m)}\,\big(\rho_{(n,m)}^{-1}(e_n\times e_m)\big)=\bigoplus_{e\in X}
\rho^{-1}(e) .
$$
\end{itemize}

\begin{rem}\label{re:IP}\textnormal{If $\rho$ has the
induction property, then $\rho$ ($=\rho_{n+m}$) and $\rho_{(n,m)}$
coincide. More precisely, if $(u,v),(u',v')\in \SG_n\times\SG_m$
such that $\rho(u\times v)=\rho(u'\times v')$, then $\rho_{(n,m)}(u
\times v) = \rho_{(n,m)}(u'\times v')$}
\end{rem}

\begin{prop}
\label{pr:FamiliesSnInduction} The graded
 maps $D$ and $\IP$ have the induction property.
\end{prop}

\begin{proof}
For  $n,m\in\NM$, let ${\mathbf c}\models n$ and ${\mathbf d}\models m$.
Let $X$ be the
subset of compositions of $n+m$ such that
\begin{equation} \label{eq:IPforD}
X_{(n,m)}\,\big(D_{(n,m)}^{-1}({\mathbf c}\times {\mathbf d})\big)\subseteq
\bigoplus_{{\mathbf e}\in X} D^{-1}({\mathbf e})
\end{equation}
and such that $X_{(n,m)}\,\big(D_{(n,m)}^{-1}({\mathbf c}\times
{\mathbf d})\big)\cap D^{-1}({\mathbf e})$ is
 non empty for all ${\mathbf e}\in X$. Since $D$ is left connected and since the $(n,m)$-components are unique, to show the equality in Eq.~(\ref{eq:IPforD}) it is
 sufficient to prove that if  $s_i\in S_{n+m}$ and
 $x(u\times v)\in
X_{(n,m)}\,\big(D_{(n,m)}^{-1}({\mathbf c}\times {\mathbf d})\big)$ are
 such that $s_ix(u\times v)\sim_D x(u\times v)$, then
 $s_i x(u\times v)\in
X_{(n,m)}\,\big(D_{(n,m)}^{-1}({\mathbf c}\times {\mathbf d})\big)$.
Using Lemma~\ref{le:Deodhar}, we reduce the problem to the case
where $s_i x = xs_j$, with $j\in [1,n-1]\cup [n+1,n+m-1]$. By
Lemma~\ref{le:Deodhar} (ii) we have
$$
1<\big|\big(x(u\times v)\big)^{-1}(i)-\big(x(u\times v)\big)^{-1}(i+1)\big| = \big|(u\times
v)^{-1}(j)-(u\times v)^{-1}(j+1)\big| .
$$
Therefore, if $j\in [1,n-1]$, we obtain $s_ju\sim_D u$, and if
$j\in [n+1,n+m-1]$, we obtain $s_{j-n}v\sim_D v$.

To prove that $\sim_\IP$ has the induction property we proceed as
above. It is sufficient to
 consider the case $i=1$ in the above proof. By definition,
$X_{(n,m)}$ is the set of permutations $x$ increasing on $[1,n]$
and on $[n+1,n+m]$. If $s_1x=x s_j$ we have  $j=1$ or $j=n+1$
 since  $x(j)=1$ and $x(j+1)=2$ by Lemma~\ref{le:Deodhar} (ii).
 In other words  $s_1 u\sim_\IP u$
 or  $s_{1}v\sim_\IP v$.
\end{proof}

The following lemma will be useful for Lemma~\ref{le:IPandCONG}.

\begin{lem}
\label{le:EquivIPIP'} If $\rho$ satisfies the Condition (IP), then
for any $\mathbf c=(c_1,\dots,c_k)\vDash n$ and
$(e_1,\dots,e_k)\in E_{c_1}\times\cdots\times E_{c_k}$,
there is $X\subseteq E_{n}$ such that
$$
X_{\mathbf c}\,\big(\rho_{\mathbf c}^{-1}(e_1\times\cdots\times e_k)\big)
=\bigoplus_{e\in X} \rho^{-1}(e) .
$$
\end{lem}

\begin{proof}  We proceed by induction on $k\geq 2$. The case $k=2$ is Condition
(IP). Assume that $k>2$. It is well known that
\begin{equation}
 \label{eq:MinimalCosetRules}
 X_\mathbf{c}=X_{(n-c_k,c_k)}\,(X_{(c_1,\dots,c_{k-1})}\times 1_{c_k}),
\end{equation}
(see for instance \cite[Lemma 2.1]{bbht}). By the induction hypothesis, there is
$X'\subseteq E_{n-c_k}$ such that
$$
X_{(c_1,\dots,c_{k-1})}\,\big(\rho_{(c_1,\dots ,c_{k-1})}^{-1}(e_1\times\cdots\times e_{k-1})\big)
 =\bigoplus_{e'\in X'} \rho^{-1}(e') .
$$
Then by Eq.~(\ref{eq:MinimalCosetRules}) we have
$$\begin{array}{rl}
X_{\mathbf c}\,\big(\rho^{-1}&(e_1\times\cdots\times e_k)\big)=
X_{(n-c_k,c_k)}\,(X_{(c_1,\dots,c_{k-1})}\times
1_{c_k})\big(\rho_\mathbf{c}^{-1}(e_1\times\cdots\times e_k)\big)\\
&= X_{(n-c_k,c_k)}\,\big(X_{(c_1,\dots,c_{k-1})}\rho_{(c_1,\dots,c_{k-1})}^{-1}
(e_1\times\cdots\times e_{k-1})\times\rho^{-1}(e_k)\big)\\
&=X_{(n-c_k,c_k)}\big(\bigoplus_{e'\in X'}
\rho_{(n-c_k,c_k)}^{-1}(e' \times e_k)\big)\\
&=\bigoplus_{e'\in X'}X_{(n-c_k,c_k)}\big(
\rho_{(n-c_k,c_k)}^{-1}(e'\times e_k)\big)\,.
\end{array}$$
Now by (IP) for each $e'\in X'$ there is $X_{e'}\subseteq E_n$ such
that $\rho_{(n-c_k,c_k)}^{-1}(e'\times e_k)=\bigoplus_{e\in X_{e'}}
\rho^{-1}(e).
$
Observe that Condition (IP) forces $X_{e_1'}\cap
X_{e_2'}=\emptyset$ if $e_1'\not=e_2'$. We set $X=\bigoplus_{e'\in
X'}X_{e'}\subseteq E_n$ and the proposition follows.
\end{proof}

\begin{lem}
\label{le:IPandCONG}%
If $\rho$ is a graded connected map with the induction property,
then for all $\sigma^\mathbf{g}\in G_n$ and $s_i\in S_n$ we have
$$
\sigma^\mathbf{g}\dot\sim^G_\rho s_i\sigma^\mathbf{g} \implies
\sigma^\mathbf{g}\dot\sim_{\rho^G} s_i \sigma^\mathbf{g}.
$$
\end{lem}

\begin{proof}
Let $\sigma^\mathbf{g}\in G_n$ and $s_i\in S_n$. Let
$\sigma^{\mathbf g}=\sigma_1^{g_1}\dots\sigma_k^{g_k}$ be the
rainbow decomposition of $\sigma^\mathbf{g}$ and let $\mathbf
c=(c_1,\dots,c_k)\vDash n$ be such that $\mathbf
r(\sigma^\mathbf{g})=\mathbf c^\mathbf{g}$. For $(u,v)\in X_\mathbf{c}\times
\SG_\mathbf{c}$, the $\mathbf c$-components of $\sigma$,
Eq.~(\ref{eq:ParabolicCompenent}) implies that
\begin{equation}
\label{eq:ParabGComponent}%
\sigma^\mathbf{g}=uv^\mathbf{g}=u
\,(\std(\sigma_1)^{g_1}\times\dots\times\std(\sigma_k)^{g_k}).
\end{equation}
Identify in the rainbow decomposition  the subword $\sigma_{\nu}$
containing $i$ as a letter and the subword $\sigma_\mu$ containing
$i+1$ as a letter. Lemma~\ref{le:Deodhar} gives us that $s_i u\in X_\mathbf{c}$ if and only if $\nu\not=\mu$. Hence, when $\nu=\mu$,  we obtain \begin{equation}
\label{eq:ParabGComponentBIS}%
s_i \sigma^\mathbf{g}= u\,(\std(\sigma_1)^{g_1}\times\dots\times
\std(s_i\sigma_\nu)^{g_\nu}\times\dots\times\std(\sigma_k)^{g_k}).
\end{equation}
Comparing Eq.~(\ref{eq:GenEquivG}) and Eq.~(\ref{eq:CongEquivG}), it is
clear that the lemma follows when $\nu\ne\mu$. For the case $\nu=\mu$ and
$\rho(s_i\sigma)=\rho(\sigma)$, since $\rho$ has the induction property,  Lemma~\ref{le:EquivIPIP'} implies
$$
\sigma,s_i\sigma\in
X_\mathbf{c}\Big(\rho_\mathbf{c}^{-1}\big(\rho_\mathbf{c}(\std(\sigma_1)\times\dots\times
\std(\sigma_\nu)\times\dots\times\std(\sigma_k)\big)\Big).
$$
Therefore comparing  Eq.~(\ref{eq:ParabGComponent}) and Eq.~(\ref{eq:ParabGComponentBIS}) we have
$\rho(\sigma_\nu)=\rho(s_i\sigma_\nu)$ and  the lemma follows in all cases.
\end{proof}

\subsection{Freeness Property}
\label{ss:MultipkicativeProperty}

 Now let us introduce the freeness property (the choice of this name will
 be explained in Corollary~\ref{co:FreeIP}).
 We say that a graded connected map $\rho\colon\SG\to\mathsf E$ or   its graded equivalence
relation $\sim_\rho$ have {\sl the freeness property} if the
following condition is satisfied:

\begin{itemize}
\item[(FP)] For any  $n,m\in\NM$, $u\in X_{(n,m)}$ and $s_j\in
S_{(n,m)}$ such that $us_j u^{-1}\in S_{n+m}$, if $v\in
\SG_{(n,m)}$ and $\rho_{(n,m)}(s_j v)=\rho_{(n,m)}(v)$,  then $\rho(us_jv)=\rho(uv)$.
\end{itemize}

\noindent
We leave it to the reader to derive the following from the  definitions.
\begin{prop}
\label{pr:FamiliesSnMP} The graded connected
 map $D$  and $\IP$ have the freeness property.
\end{prop}

The following lemma will be useful to show the converse to Lemma~\ref{le:IPandCONG}.

\begin{lem}
\label{le:EquivMPMP'} If $\rho$ satisfies the freeness property, then for any
$\mathbf c=(c_1,\dots,c_k)\vDash n$, $u\in X_{(n,m)}$ and
$s_j\in S_{\mathbf c}$ such that $us_j u^{-1}\in S_{n}$, we have that $v\in
\SG_{\mathbf c}$ and $\rho_\mathsf{c}(s_j v)=\rho_\mathsf{c}(v)$
implies
$\rho(us_jv)=\rho(uv)$.
\end{lem}

\begin{proof}  Let  $v\in \SG_{\mathbf c}$ be such that $\rho(s_j
v)=\rho(v)$ and write $s_i = us_j u^{-1}\in S_n$. We prove that
$\rho(us_jv)=\rho(uv)$ by induction on $k$. The case $k=2$ is
(FP). Assume that $k>2$ and write $v=v_1\times\dots\times v_k$.
Using Eq.~(\ref{eq:MinimalCosetRules}) we write
 $u=u'(u''\times 1_{c_k})$ with $u'\in X_{(n-c_{k},c_k)}$ and $u''\in
 X_{(c_1,\dots,c_{k-1})}$. Since $u'\in X_{(n-c_k,c_k)}$, Lemma~\ref{le:Deodhar} gives us that $s_i u'\in
 X_{(n-c_k,c_k)}$ or there is $s_p\in S_{(n-c_k,c_k)}$ such that
 $s_i u'=u' s_p$. The second case forces $s_p (u''\times 1_{c_k}) = (u''\times
 1_{c_k}) s_j$. Since we have $\ell(s_i u)=\ell(us_j)=\ell(u)+1$ and
 $\ell\big(u'(u''\times 1_{c_k})\big)=\ell(u')+\ell(u''\times 1_{c_k})$,
 the first case would  imply that
 $\ell(s_iu')=\ell(u')+1$. This would imply
$$
\ell(u)=\ell(s_i u s_j)=\ell(s_iu' (u''\times 1_{c_k})s_j)
=\ell(s_iu')+\ell (u''\times 1_{c_k})+1>\ell(u),
$$
which is a contradiction. Hence $s_p= (u''\times 1_{c_k})s_j(u''\times 1_{c_k})^{-1}\in
S_{(n-c_k,c_k)}$.
For $s_j\in\SG_{n-c_k}$ the induction hypothesis gives us $\rho\big(u'' s_j (v_1\times \dots\times
 v_{k-1})\big)=\rho\big(u'' (v_1\times \dots\times
 v_{k-1})\big)$. The case when $s_{j}\in \SG_{(1,\dots,1,c_k)}$ is
 trivial. Then
 $\rho_\mathbf{c}( (u''\times 1_{c_k}) s_j v)= \rho_\mathbf{c}( (u''\times 1_{c_k})
 v)$. The result follows  using (FP) and the fact that $(u''\times 1_{c_k}) s_j=s_p (u''\times
 1_{c_k})$ and $s_iu'=u's_p$.
\end{proof}

\begin{lem}
\label{le:MPandCONG}%
If $\rho$ is a graded connected map satisfying Condition (FP),
then for all $\sigma^\mathbf{g}\in G_n$ and $s_i\in S_n$ we have
$$
\sigma^\mathbf{g} \dot\sim_{\rho^G} s_i\sigma^\mathbf{g} \implies
\sigma^\mathbf{g}\dot\sim^G_\rho s_i \sigma^\mathbf{g}.
$$
\end{lem}

\begin{proof} We use the same notation as in the proof of
Lemma~\ref{le:IPandCONG}. Again, comparing
Eq.~(\ref{eq:GenEquivG}) and Eq.~(\ref{eq:CongEquivG}), it is
clear that we only have to consider the case $\nu=\mu$.
Lemma~\ref{le:Deodhar} implies that $s_i=us_ju^{-1}$ with $s_j\in
S_\mathbf{c}$. Comparing Eq.~(\ref{eq:ParabGComponent}) and
Eq.~(\ref{eq:ParabGComponentBIS}), we obtain that
$\std(s_i\sigma_\nu)=s_{j_0}\std(\sigma_\nu)$ with
$j_0=j-(c_1+\dots+c_\nu)$. Then the condition
$\rho(s_i\sigma_\nu)=\rho(\sigma_\nu)$ is equivalent to
$\rho(s_{j_0}\std(\sigma_\nu))=\rho(\sigma_\nu)$. That is
$\rho_{\mathbf c}(s_jv)=\rho_\mathbf{c}(v)$ and the result follows
from Lemma~\ref{le:EquivMPMP'}.
\end{proof}

\begin{prop}
\label{pr:FPandIP}%
 If $\rho$ is a graded connected map satisfying
Conditions (IP) and (FP), then $\sim_{\rho^G}$ and
$\dot\sim^G_\rho$ are equal. In particular, the
graded connected maps $\rho_G$ and $\rho^G$ induce the same equivalence relations.
\end{prop}

For instance, consider the  graded connected descent map $D_G$
induced  from the descent (surjective) map $D\colon\SG\to\Comp$.
Applying Propositions~\ref{pr:AtkinsonDescent} and
~\ref{pr:FPandIP}, we get

\begin{cor}
\label{co:GDescentEquiv} For any $\alpha,\beta\in G_n$,
$D_G(\alpha)=D_G(\beta)$ if and only if $\alpha\sim_D^G \beta$.
\end{cor}

The description of $D_G$ given by the elementary relation
$\dot\sim^G_D$ is the principal ingredient of the proof given in
\cite{baumann hohlweg} that $\Sigma_n(G)$ is a subalgebra of $\ZM
G_n$.

\subsection{Coloured peak algebras}
\label{ss:ColouredPeakAlg}

Let $\sigma^\mathbf{g}=\sigma_1^{g_1}\dots \sigma_k^{g_k}$ be the
rainbow decomposition of $\sigma^\mathbf{g}\in G_n$. The {\sl
$G$-peak composition of $\sigma^\mathbf{g}$} is
$$
\IP_G(\sigma^\mathbf{g})= \big( \IP(\sigma_1)^{g_1},\dots,
\IP(\sigma_k)^{g_k} \big).
$$
Denote by $\IPi_n(G)$ the set of $G$-peak compositions of $G_n$.
Then $(\IPi^G)^\ast=\bigoplus_{n\in\NM} \IPi_n(G)$.

In the group algebra $\ZM G_n$, we form  the elements
$
 \ip^G_{\mathbf
c^{\mathbf g}} = \sum_{\IP_G(\alpha)=\mathbf c^\mathbf{g}} \alpha,
$
where $\mathbf c^\mathbf{g}\in\IPi_n(G)$. We then consider the
subspace
$$
\IPeak_n(G)=\bigoplus_{\mathbf c^{\mathbf g}\in\IPi_n(G)}\ZM
\ip^G_{\mathbf c^{\mathbf g}}  .
$$

\begin{thm}
\label{th:ColPeakAlgebras} If $G$ is a finite abelian group and $n\in\NM$,
 then  $\IPeak_n(G)$ is a (non-unitary) subalgebra of $\Peak_n(G)$ called
 the $G$\textup{-peak algebra}.
\end{thm}

We proceed as in \S~\ref{ss:DescentSets} and the theorem will follow from the next two lemmas.
First, Propositions~\ref{pr:FPandIP} and
 Corollary~\ref{co:GDescentEquiv}, combined with Eq.~(\ref{eq:InteriorPeaksInDescent}), show that

\begin{lem}
\label{le:GPeaks} Let $\alpha,\beta\in G_n$.
\begin{enumerate}

\item[(i)] $\alpha\sim_D^G \beta \implies \alpha\sim_\IP^G \beta$;

\item[(ii)] $\IP_G(\alpha)=\IP_G(\beta) \iff \alpha\sim_\IP^G
\beta$.
\end{enumerate}
\end{lem}

For $\alpha\in G_n$ and $\mathbf c^\mathbf{g},\mathbf
d^\mathbf{h}\in\IPi_n(G)$, we set
$$
\IP_G(\mathbf c^\mathbf{g},\mathbf d^\mathbf{h},\alpha) =
\set{(\sigma^\mathbf{g},\tau^\mathbf{h}) \in G_n \times G_n \st
\IP_G(\sigma^\mathbf{g})=\mathbf c^\mathbf{g},\
\IP_G(\tau^\mathbf{h})=\mathbf d^\mathbf{h},\
\sigma^\mathbf{g}\tau^\mathbf{h}=\alpha }.
$$
 Fix $s_i \in S_n$ and  let
$(\sigma^\mathbf{g},\tau^\mathbf{h}) \in G_n \times G_n$. Then
$$
\psi_i^{\IP_G} (\sigma^\mathbf{g},\tau^\mathbf{h}) =
\left\{\begin{array}{ll}
(s_i \sigma^\mathbf{g},\tau^\mathbf{h}) &\textrm{if } s_i \sigma^\mathbf{g}\dot\sim_{\IP_G} \sigma^\mathbf{g}\\ \\
(\sigma^\mathbf{g},\sigma^{-1}s_i\sigma\tau^\mathbf{h})&\textrm{otherwise}
\end{array}\right.
$$
is an involution on $G_n \times G_n$.

\begin{lem}
\label{le:PsiGPeak} If $\alpha \in G_n$ and $s_i \in S_n$ are
such that $\alpha \dot\sim_\IP^G s_i \alpha$, then
$$
\psi^{\IP_G}_i\big(\IP_G(\mathbf c^\mathbf{g},\mathbf
d^\mathbf{h},\alpha)\big)=\IP_G(\mathbf c^\mathbf{g},\mathbf
d^\mathbf{h},s_i\alpha),\quad\textrm{for all }\mathbf
c^\mathbf{g},\mathbf d^\mathbf{h} \in \IPi_n(G).
$$
\end{lem}
\begin{proof} We can assume that $i$ and $i+1$ are in the same subfactor
of the rainbow decomposition of $\sigma^\mathbf{g}$ and that
$s_i\sigma \not\sim_\IP \sigma$. That is, there is $s_j\in S_n$
such that $s_i\sigma=\sigma s_j$ and $s_j$ stabilize $\mathbf g$.
Therefore $s_i\sigma^\mathbf{g} = \sigma^\mathbf{g} s_j$. If $j$
and $j+1$ are not in the same subfactor of the rainbow
decomposition of $\tau^\mathbf{h}$ then
$s_j\tau^\mathbf{h}\dot\sim_{\IP_G} \tau^\mathbf{h}$ by
definition. Otherwise, observe that $|\alpha|\dot\sim_\IP s_i
|\alpha|$ and
$$
\psi_i^{\IP_G}
(\sigma^\mathbf{g},\tau^\mathbf{h})=\psi_i^\IP(\sigma,\tau)^{(\mathbf
g,\mathbf h)}
$$
where $(u,v)^{(\mathbf g,\mathbf h)}=(u^\mathbf{g},v^\mathbf{h})$.
We conclude that $s_j\tau^{\mathbf h}\dot\sim_{\IP_G}
\tau^\mathbf{h}$ since  $\Psi_i^\IP (\sigma,\tau) \in$\break
$\IP(\IP(\sigma),\IP(\tau),s_i|\alpha|)$
 (Lemma~\ref{le:PsiIntPeak}).
\end{proof}

%
%

\section{graded maps and bialgebra structures}
\label{se:HopfAlgebras}

Henceforth we consider a field $\k$ of
characteristic $0$.

\subsection{Hopf algebra of coloured permutations}
\label{ss:HopfColouredPermutations}

 Let $\rho\colon\mathsf G\to\mathsf E$ be a graded map and $n\in\NM$.
We denote by $\mathcal E(\rho)_n$ the subspace of the group
algebra $\k G_n$ spanned by
$
\mathsf \b_e^{\rho}=\sum_{\rho(\alpha)=e} \alpha$ for each $e
\in E_n.
$
Then
$$
\mathcal E(\rho)=\bigoplus_{n\geq 0} \mathcal E(\rho)_n
$$
is a graded vector subspace of the graded vector space
$
\mathcal \k[\mathsf G]=\mathcal
E(\textrm{Id}_\mathsf{G})=\bigoplus_{n\geq 0} \k G_n .
$

 We  recall the Hopf algebra structure on $ \k[\mathsf G]$
given in \cite{baumann hohlweg}.
The {\sl product} $*\colon\k G_n\otimes \k G_m \to \k G_{n+m}$ is
given by
$$
\alpha * \beta = x_{n,m}(\alpha \times \beta) \in \k G_{n+m}.
$$

Let $\alpha \in G_n$. For each $i \in [0,n]$, we denote
by $\alpha_{(i)}\times\alpha_{(n-i)}$ the unique element of
$G_{i,n-i}$ such that $\Big(u, (\alpha_{(i)}\times\alpha_{(n-i)})^{-1}\Big)$ are
the $(i,n-i)$-components of $\alpha^{-1}$. The {\sl coproduct}
$\Delta\colon\k[\mathsf G]\to\k[\mathsf G]\otimes \k[\mathsf G]$ is
given by
$$
\Delta(\alpha) = \sum_{i=0}^n \alpha_{(i)} \otimes_\ZM
\alpha_{(n-i)} \in \k[\mathsf G]\otimes \k[\mathsf G].
$$

For instance, for $g\not=h$ in $G$, $$ 1^g 2^h * 2^h 1^g =1^g 2^h 4^h 3^g+ 1^g 3^h 4^h
2^g+ 1^g4^h3^h 2^g+ 2^g 3^h 4^h 1^g + 2^g 4^h 3^h 1^g+3^g 4^h 2^h 1^g;$$
$$
 \Delta (2^g 3^h 1^h 4^g)= \delta\otimes 2^g 3^h 1^h
4^g+ 1^h \otimes 1^g 2^h  3^g + 2^g 1^h\otimes 1^h 2^g +
2^g 3^h 1^h \otimes  1^g+ 2^g 3^h1^h 4^g \otimes\delta .
$$

In \cite{baumann hohlweg}, the authors prove that $(\k[\mathsf
G],*,\Delta)$ is a graded connected Hopf algebra. Moreover, the
assignment $G\rightsquigarrow \k[\mathsf G]$ is a covariant
functor from the category of finite abelian groups to the category
of graded Hopf algebras. The case $G=\set{1}$ is due to Malvenuto and
Reutenauer \cite{malvenuto reutenauer}. Considering a finite
abelian group $G$ and the morphism $G\to\set{1}$, at the
level of graded Hopf algebras, we get  the forgetful (surjective) map
$f\colon\k[\mathsf G]\to \k[\SG]$ (obtained by extending linearly to
$\k[\mathsf G]$ the absolute value map $|\,.\,|$). We also have an injective morphism induced from the injective map $\{1\}\hookrightarrow G$.

We construct more Hopf (sub)algebras using our general theory of graded maps on $\mathsf G$. For example:
\begin{enumerate}

\item The {\sl Solomon descent algebra} $ \Sigma=\mathcal E
(D\colon\SG\to (\Comp,\|\,.\,\|))$ is a Hopf subalgebra of $\k[\SG]$
\cite{malvenuto reutenauer}.

\item The Mantaci-Reutenauer algebra  $ \Sigma(G)=\mathcal E(D_G)$
is a Hopf subalgebra of $\k[\mathsf G]$. Moreover, the assignment
$G\rightsquigarrow \Sigma(G)$ is also a contravariant functor from
the category of finite abelian groups to the category of graded
Hopf algebras \cite{baumann hohlweg}.

\item The {\sl Peak algebra} $\IPeak=\mathcal E(\IP\colon\SG\to
(\Comp,\|\,.\,\|))$ is a Hopf subalgebra of $\k[\SG]$
\cite{aguiar bergeron nyman, bht}.

\item The {\sl planar binary tree $\UT(\sigma)$ of $\sigma\in
\SG_n$} is defined as
 follows. A given $\sigma\in \SG_n$ can
 be viewed as the concatenation of two injective words
  separated by the smallest letter $1$:
  $\sigma=w' \cdot 1 \cdot w''$. We construct $\UT$ on
  permutation by induction: $\UT(\sigma)=\UT(\std(w'))\vee
  \UT(\std(w''))$.  Here $T'\vee T''$ means {\sl the grafting} of $T'$ and $T''$,
see \cite{loday ronco}. Denoting $Y_n$ the set of planar binary
trees with $n$ vertices,
 we get a surjective graded map $\UT\colon\SG\to\mathsf Y$, where
 $\mathsf Y=\bigoplus_{n\in\NM} Y_n$. Then $\mathcal E(\UT)$ is the Hopf subalgebra of $\k[\SG]$,
 called the {\sl Loday-Ronco algebra}  \cite{loday ronco}.

\end{enumerate}
We see that the example (2) is a $G$-coloured version of (1). In Section~\ref{se:Appli} we will present a coloured version of (3) and (4). For this, we study which properties of $\rho$ are required so that $\mathcal E(\rho)$  is a Hopf subalgebra.

\subsection{Induction on equivalence classes}
\label{ss:IP}

Let $\rho\colon\mathsf G\to\mathsf E$ be a graded map. The
definition of the induction property (IP) in
\S\ref{ss:InductionProperty} can be applied to $\rho$ viewing the product in $\mathsf G$ instead of in $\SG$. Since the map $X_{n,m}\times G_{n,m}\to
G_{n+m}$, sending $(u,v^\mathbf{g})$ to $uv^{\mathbf g}$, is a
bijection, we get the following:

\begin{prop}
\label{pr:HopfInduction}
A graded map $\rho\colon\mathsf G\to\mathsf E$
satisfies Condition (IP)
 if and only if for all $n,m\in\NM$ and for all
$(e_1,e_2)\in E_n\times E_m$, there is $X\subseteq E_{n+m}$ such
that
$$
\b^{\rho}_{e_1} * \b^{\rho}_{e_2} = \sum_{e\in X}\b^{\rho}_e .
$$
\end{prop}

\begin{cor}
\label{co:HopfInduction}
 Let $\rho$ be a graded map on $\mathsf G$.  If $\rho$
has the induction property (IP), then $\mathcal E (\rho)$ is a
subalgebra of $(\k[\mathsf G],*)$.
\end{cor}

\subsection{Restriction on equivalence classes}
\label{ss:RP}

The graded map $\rho\colon\mathsf G\to \mathsf E$ has {\sl the
restriction property} if it satisfies  the following condition:

\begin{itemize}

\item[(RP)] For any $n,m\in\NM$ and  $u\in X_{n,m}$, if
$\alpha_1,\alpha_2\in G_n$ and $\beta_1,\beta_2\in G_m$ are such that
$\rho(\alpha_1)=\rho(\alpha_2)$ and  $\rho(\beta_1)=\rho(\beta_2)$, then $$\rho\big(
(\alpha_1\times \beta_1)u^{-1}\big) = \rho\big((\alpha_2\times
\beta_2)u^{-1}\big).$$
\end{itemize}

\begin{rem}\label{re:RP} \textnormal{If $\rho$ satisfies (RP), we obviously
have:}

\textnormal{(i) $\rho(\alpha_1\times \beta_1)= \rho(\alpha_2\times
\beta_2)$, since $1_{n+m}\in X_{n,m}$. But we do  not have that
$\rho_{(n,m)}(\alpha_1\times \beta_1)= \rho_{(n,m)}(\alpha_2\times
\beta_2)$, this  comes from the induction property
(see Remark~\ref{re:IP}).}

\textnormal{(ii) For $n\in \NM$,  $\alpha,\beta\in G_n$  and $i\in
[0,n]$, if
 $\rho(\alpha_{(i)})= \rho(\beta_{(i)})$  and
$\rho(\alpha_{(n-i)})=\rho(\beta_{(n-i)})$, then
$\rho(\alpha_{(i)}\times\alpha_{(n-i)})=\rho(\beta_{(i)}\times\beta_{(n-i)}).$}
\end{rem}

\begin{lem}
\label{le:EquivHaveRP}
 The graded connected maps $D$ and  $\IP$ satisfy Condition (RP).
\end{lem}
\begin{proof}
Let $s_i\in S_{(n,m)}$ and $\sigma\times \tau\in \SG_{(n,m)}$ be such
that $s_i\sigma \dot\sim_D \sigma$ when $s_i\in S_n$ and
$s_{i-n}\tau\dot\sim_D \tau$ when $s_{i-n}\in S_m$. In both cases
 we have $|(\sigma\times \tau)^{-1}(i+1)-(\sigma\times
\tau)^{-1}(i)|>1$. For $u\in X_{(n,m)}$ we have
$$
|((\sigma\times \tau)u^{-1})^{-1}(i+1)-((\sigma\times
\tau)u^{-1})^{-1}(i)|=|u(\sigma\times
\tau)^{-1}(i+1)-u(\sigma\times \tau)^{-1}(i)|
$$
Observing that $(\sigma\times \tau)^{-1}(i+1)$ and
$(\sigma\times \tau)^{-1}(i)$ are elements of either  $[1,n]$ or
$[n+1,n+m]$, we have that Eq~(\ref{eq:EquivDescentSn}) holds.
The lemma follows from the fact that $D$ is left connected.

For $\IP$ the proof is similar.
\end{proof}

We now show that if $\rho$ has the restriction property, then
$\mathcal E (\rho)$ is a subcoalgebra of $(\k[\mathsf G],\Delta)$.
For $n\in \NM$, $e\in E_n$,  $i\in[0,n]$ and
$(\beta_1,\beta_2)\in E_i\times E_{n-i}$ we set
\begin{equation}
\label{eq:CoeffRP}%
 A^e_{\beta_1,\beta_2} =\set{\alpha\in
\rho^{-1}(e)\st \ \alpha_{(i)}\times \alpha_{(n-i)} =\beta_1\times
\beta_2}.
\end{equation}

\begin{lem}
\label{le:RP}%
 Let $\rho$ be a graded map on $\mathsf G$ satisfying (RP).
 For $n\in \NM$,  $e\in E_n$ and $i\in
[0,n]$, if $\beta_1',\beta_1''\in E_i$ and $\beta_2',\beta_2''\in
E_{n-i}$ are such that
$
\rho(\beta_1')=\rho(\beta_1'') \textrm{ and } \rho(\beta_2')=
\rho(\beta_2''),
$
then $A^e_{\beta_1',\beta_2'}$ and $A^e_{\beta_1'',\beta_2''}$ are
in bijection.
\end{lem}

\begin{proof} If $\alpha'\in A^e_{\beta_1',\beta_2'}$,
then $\alpha'=(\beta_1'\times \beta_2') u^{-1}$ with $u\in
X_{(i,n-i)}$. Since the $(i,n-i)$-components are unique, we
define a map from $A^e_{\beta_1',\beta_2'}$ to $G_n$ by
sending $\alpha'$ to $\alpha''= (\beta''_1\times \beta''_2)
u^{-1}$. By the restriction property, Eq.~(\ref{eq:CoeffRP})  and
Remark~\ref{re:RP}, we have
$$
\rho(\alpha'')=\rho( (\beta_1''\times \beta_2'') u^{-1} ) =\rho(
(\beta_1'\times \beta_2') u^{-1} ) =\rho(\alpha')=e.
$$
In other words, $\alpha''\in A^e_{\beta_1'',\beta_2''}$. By
uniqueness of the $(i,n-i)$-components again, and by symmetry of
the construction, we get a bijection. If $A^e_{\beta_1',\beta_2'}$
is empty, the above discussion shows that
$A^e_{\beta_1'',\beta_2''}$ is also empty.
\end{proof}

Let $e\in E_n$. Lemma~\ref{le:RP} allows us to define for  all
$i\in [0,n]$ and for all $(e_1,e_2)\in E_i \times E_{n-i}$ the
integer
$$
a^e_{e_1,e_2}=|A^e_{\beta_1,\beta_2}|,\quad\textrm{for }
(\beta_1,\beta_2) \in \rho^{-1}(e_1)\times\rho^{-1}(e_2).
$$

\begin{thm}
\label{th:RP}
Let $\rho$ be a graded map on $\mathsf G$ satisfying (RP) and $n\in\NM$.
For all $e\in E_n$ we have
$$
\Delta(\b^{\rho}_e)=\sum_{i=0}^n \sum_{(e_1,e_2)\in E_i\times
E_{n-i}} a^e_{e_1,e_2}\, \b^{\rho}_{e_1}\otimes\b^{\rho}_{e_2} .
$$
 In particular,  $\mathcal E(\rho)$ is a
subcoalgebra of $(\k[\mathsf G],\Delta)$.
\end{thm}

\begin{proof} For each $i\in [0,n]$, we have
$$
\rho^{-1}(e)=\bigoplus_{(\beta_1,\beta_2)\in G_i\times G_{n-i}}
A^e_{\beta_1,\beta_2}.
$$
Then
\begin{eqnarray*}
\Delta(\b^{\rho}_e)
      &=& \sum_{i=0}^n
      \sum_{\alpha\in \rho^{-1}(e)}\alpha_{(i)}\otimes
      \alpha_{(n-i)}\\
     &=& \sum_{i=0}^n \sum_{(\beta_1,\beta_2)\in G_i\times G_{n-i}}
      \sum_{\alpha\in A^e_{\beta_1,\beta_2}}\beta_1\otimes
      \beta_2\\
     &=&\sum_{i=0}^n \sum_{(e_1,e_2)\in E_i\times E_{n-i}}
        \sum_{(\beta_1,\beta_2)\in \rho^{-1}(e_1)\times\rho^{-1}(e_2)}
      a^e_{e_1,e_2}\,\beta_1\otimes
      \beta_2\\
     &=&\sum_{i=0}^n \sum_{(e_1,e_2)\in E_i\times E_{n-i}} a^e_{e_1,e_2}\,
  \b^{\rho}_{e_1}\otimes\b^{\rho}_{e_2}.
\end{eqnarray*}
\end{proof}

\subsection{Generated connected graded maps and bialgebra structures}
\label{ss:HopfGenerated}

Recall the definition of $\rho_G$ in \S\ref{ss:GenEquivG}. The following
theorem gives an automatic way to build coloured Hopf algebras.

\begin{thm}\label{th:Generating} Let $G$ be a finite abelian group and
let $\rho\colon\SG\to \mathsf E$ be a connected graded map.

\begin{enumerate}

\item[(i)] If $\rho$ has the induction property, then $\rho_G$ has
the induction property.

\item[(ii)] If $\rho$ has the restriction property, then $\rho_G$
has the restriction property.

\end{enumerate}
\end{thm}

\begin{proof}  (i)  Let  $n,m\in\NM$ and let $(e_n(G),e_m(G))
\in E_n(G)\times E_m(G)$. Let $X\subseteq E_{n+m}(G)$ be the
subset  such that
\begin{equation}\label{eq:IPforGen}
X_{(n,m)}\,\big((\rho_G)_{(n,m)}^{-1}(e_n(G)\times e_m(G)\big)\subseteq
\bigoplus_{e(G)\in X} \rho_G^{-1}(e(G))
\end{equation}
and such that $X_{(n,m)}\,\big((\rho_G)_{(n,m)}^{-1}(e_n(G)\times
e_m(G)\big)\cap\rho_G^{-1}(e(G))$ is
 non empty for all $e(G)\in X$. Using Eq.~(\ref{eq:CongEquivG}), since $\rho_G$ is connected, to show the equality in Eq.~(\ref{eq:IPforGen})   it is
 sufficient to prove that if  $s_i\in S_{n+m}$ and $\sigma^\mathbf{g}\in
X_{(n,m)}\,\big((\rho_G)_{(n,m)}^{-1}(e_n(G)\times e_m(G)\big)$ are
 such that $s_i\sigma^\mathbf{g}\dot\sim_\rho^G \sigma^\mathbf{g}$, then
 $s_i\sigma^\mathbf{g}\in
X_{(n,m)}\big((\rho_G)_{(n,m)}^{-1}(e_n(G)\times e_m(G)\big)$. Let
$(u,\tau_1 \times \tau_2)\in X_{(n,m)}\times \SG_{(n,m)}$ be the
$(n,m)$-components of $\sigma$. Then
$\sigma^\mathbf{g}=u(\tau_1 \times \tau_2)^{\mathbf g_1\cdot \mathbf  g_2}$.
 Using Lemma~\ref{le:Deodhar} we obtain either $s_i u\in X_{(n,m)}$
and the proof is done, or there is $j\in S_{(n,m)}$ such that $s_iu=us_j$
 (and $u(j)=i$ and $u(j+1)=i+1$).
 Assume first that $s_j\in S_n$.
  If $j$ and $j+1$ are not in the same subfactor of the rainbow decomposition
 of $\tau_1^{\mathbf g_1}$ the proof is done. If $j$ and $j+1$ are in
the same subfactor then $i$ and $i+1$ are in the same subfactor of
 the rainbow decomposition of $\sigma^\mathbf{g}$ since
$u(j)=i$ and $u(j+1)=i+1$. Then by definition $\rho(s_i\sigma) =\rho(\sigma)$.
 The fact that $\rho$ has the induction property forces
$\rho(s_j\tau_1)=\rho(\tau_1)$. Hence
 $s_j\tau_1^{\mathbf g_1}\dot\sim_\rho^G \tau_1^{\mathbf g_1}$ and
 the proof follows. Proceed similarly
 if $j-n \in S_m$.

(ii) Let $s_i\in S_{(n,m)}$, $\sigma^\mathbf{g}\times
\tau^\mathbf{h}\in G_n\times G_m$ and $u\in X_{(n,m)}$. Observing
that $i$ and $i+1$ are in the same subfactor of the rainbow
decomposition of $\sigma^\mathbf{g}\times \tau^\mathbf{h}$ if and
only if $i$ and $i+1$ are in the same subfactor of the rainbow
decomposition of $(\sigma^\mathbf{g}\times
\tau^\mathbf{h})u^{-1}$, this case follows from definitions.
\end{proof}

Let $\rho\colon\SG\to \mathsf E$ be a induced graded map satisfying
(FP) and (IP). We have a nice description of the
equivalence classes
 of $\rho_G$ by Proposition~\ref{pr:FPandIP}. For $g\in G$, using the map in
Eq.~(\ref{eq:BijInducedEquiv}), we define a monomorphism of graded
vector spaces
\begin{eqnarray*}
\mu_g\colon\, \mathcal E(\rho) &\to& \mathcal E(\rho_G)\\
\b^\rho_{e} &\mapsto&\b^{\rho_G}_{e^g}.
\end{eqnarray*}

\begin{cor}
\label{co:FreeIP}%
Let $\rho\colon\SG\to \mathsf E$ be a  graded connected  map satisfying
(FP) and (IP). If  $\mathcal
E(\rho)$ is freely generated by $M$, then $\mathcal E(\rho_G)$ is
a subalgebra of $\k[\mathsf G]$ freely generated by
$$
M(G)=\bigoplus_{g\in G} \mu_g \left( M \right).
$$
\end{cor}

\begin{proof} Observe that for $g_1,g_2 \in G$ and $e_1 , e_2 \in \mathsf E$
 we have
$$
\b^{\rho_G}_{e_1^{g_1}}* \b^{\rho_G}_{e_2^{g_2}} = \left\{\begin{array}{ll} %
\b^{\rho_G}_{e_1^{g_1}\cdot e_2^{g_2}}&\textrm{if } g_1\not= g_2\,,\\ \\
\mu_{g_1}\big(\b_{e_1}^{\rho} * \b^{\rho}_{e_2}\big)& \textrm{if } g_1= g_2\,.
\end{array}\right.
$$
Conclude by induction.
\end{proof}

\begin{cor}
\label{co:Fonctoriality}%
Let $\rho$ be a graded connected map satisfying  (FP) and (IP). If
$\mathcal E(\rho)$ is free, then the assignment $G\rightsquigarrow
\mathcal E(\rho_G)$ is a covariant functor from the category of
finite abelian groups to the category of free graded algebras.
\end{cor}

\begin{proof} Recall that  the
assignment $G\rightsquigarrow \k[\mathsf G]$ is a covariant
functor from the category of finite abelian groups to the category
of graded Hopf algebras \cite{baumann hohlweg}. More precisely if
$f\colon G\to G'$ is a group homomorphism, we define a Hopf algebra
homomorphism $f_*$ by sending $\sigma^\mathbf{g}\in\k[\mathsf G]$
to
 $\sigma^{f(\mathbf{g})}\in\k[\mathsf G']$
(where $f((g_1,\dots,g_n))=(f(g_1),\dots,f(g_n))$).
Then we have $f_*\circ\mu_g=\mu_{f(g)}$ for all $g\in G$. In other
words, $f_* (M(G))\subseteq M(G')$. The corollary follows.
\end{proof}

%
%

\section{Applications}
\label{se:Appli}

\subsection{Coloured Bialgebras of  peaks and trees}
\label{ss:ColouredBialgPeakTree}

A first application of the above theory is the following theorem:

\begin{thm}
\label{th:PeakHopf} %
Let $G$ be an abelian group.
\begin{enumerate}

\item[(i)] The graded space
$$
\IPeak(G)=\bigoplus_{n\in\NM} \IPeak_n(G)
$$
is a subalgebra of $(\Sigma(G),*)$ freely generated by
$$
\set{\ip^{G}_{(n)^g}\st g\in G,\ n\textrm{ odd}}.
$$

\item[(ii)] The graded algebra $\IPeak(G)$ is a Hopf subalgebra of
$(\Sigma(G),*,\Delta)$. Moreover, the assignment
$G\rightsquigarrow \IPeak(G)$ is a covariant functor from the
category of finite abelian groups to the category of graded
connected Hopf algebras.
\end{enumerate}
\end{thm}
\begin{proof} The fact that $\IPeak(G)$ is a Hopf subalgebra
of $(\Sigma(G),*,\Delta)$ comes from Proposition~\ref{pr:FamiliesSnInduction}
and \ref{le:EquivHaveRP},
Theorems~\ref{th:RP} and \ref{th:Generating}, and Corollary~\ref{co:HopfInduction}.
The peak algebra $\IPeak$ is freely generated by $\ip_{(n)}$ for $n$ odd
 \cite[Th. 5.4]{bmws}. Then the freeness  follows from Corollary~\ref{co:FreeIP}, and the functorial
property follows from Corollary~\ref{co:Fonctoriality}.
\end{proof}

We next consider a $G$-colouring of the Loday-Ronco Hopf algebra of trees \cite{loday ronco}.
Let $\sigma^\mathbf{g}=\sigma_1^{g_1}\dots \sigma_k^{g_k}$ be the
rainbow decomposition of $\sigma^\mathbf{g}\in G_n$. The {\sl
$G$-sequence of trees of $\sigma^\mathbf{g}$} is
$$
\mathcal T_G(\sigma^\mathbf{g})= \big( \mathcal
T(\sigma_1)^{g_1},\dots, \mathcal T(\sigma_k)^{g_k} \big).
$$

Hivert, Novelli and Thibon \cite{hivert novelli thibon,hivert
novelli thibon 1,hivert novelli thibon 2} have shown that the
graded map $\mathcal T$ is connected. The {\sl
dual-sylvester equivalence} is defined by
$$
s_i \sigma \dotsimsylv \sigma \quad\iff\quad \exists\, k\in\,
[w^{-1}(i),w^{-1}(i+1)], \ \sigma(k)>i+1
$$
They  also show that $\UT(w)=\UT(g)\iff w\simsylv g$. The
following lemma has been indirectly proved in \cite{hivert novelli
thibon}.

\begin{lem}
\label{le:SylvIPRP}%
The graded connected map $\mathcal T$ satisfies (MP), (IP) and
(RP).
\end{lem}
\begin{proof} Proceed as  in Proposition~\ref{pr:FamiliesSnInduction},
\ref{pr:FamiliesSnMP} and \ref{le:EquivHaveRP}.
\end{proof}

\begin{thm}
\label{th:ColouredTrees}%
The graded space $\mathcal E\big(\mathcal T_G\big)$ is a graded
connected Hopf subalgebra of $\k[\mathsf G]$ containing
$\Sigma(G)$. As an algebra, it is freely generated by
$$
 \set{\b_{(|\vee T)^g}^{\mathcal T_G}\st g\in G,\ T\in Y_{n-1}},
$$
where $|\in Y_1$ is the unique tree of degree 1 and $\vee$ is the
grafting operation (see \cite{loday ronco}). Moreover, the
assignment $G\rightsquigarrow \mathcal E\big(\mathcal T_G\big)$ is
a covariant functor from the category of finite abelian groups to
the category of graded connected Hopf algebras.
\end{thm}
\begin{proof} The fact that $\mathcal E\big(\mathcal T_G\big)$ is a Hopf subalgebra
of $(\k[\mathsf G],*,\Delta)$ comes from Lem\-ma~\ref{le:SylvIPRP},
Theorems~\ref{th:RP} and \ref{th:Generating}, and Corollary~\ref{co:HopfInduction}.
The Loday-Ronco algebra of trees $\mathcal E(\mathcal T)$ is freely generated by
$\b_{|\vee T}^{\mathcal T}$ for  $T\in Y_{n-1}$
  \cite[Theorem 3.8]{loday ronco}. The freeness  follows from Corollary~\ref{co:FreeIP}, and the functorial
property follows from Corollary~\ref{co:Fonctoriality}.
\end{proof}

\begin{rem}\textnormal{
We could also apply this theory to the Knuth relations and the
Hopf algebra of tableaux \cite{poirier reutenauer}. This would
lead to a different $G$-coloured version than the one defined in
\cite[\S5.5]{baumann hohlweg}.}
\end{rem}

\subsection{Exterior peaks in the symmetric group}
\label{ss:PeakSets}

Let $\sigma\in\SG_n$, the {\sl set of exterior peaks of}
 $\sigma$ is
$$
\peak(\sigma)=\set{i\in [1,n-1]\st
\sigma(i-1)<\sigma(i)>\sigma(i+1)}
$$
where we set $\sigma(0)=0$ \cite[Definition 3.1]{aguiar bergeron
nyman}.
  Notice that $1\in \peak(\sigma)$ if and
only if $1\in \Des(\sigma)$. Similarly to peak sets, all the
subsets of $[1,n-1]$ are not exterior peak sets. In fact $I\subseteq
[1,n-1]$ is an exterior peak set if and only if $I$ satisfies the
condition: if $i\in I\setminus\set{1}$, then $i-1\notin I$. The
{\sl peak composition $P(\sigma)$ of $\sigma\in\SG_n$} is defined
by Eq.~(\ref{eq:BijCompSubset2}): $P(\sigma)=\mathbf
c_{\peak(\sigma)}$. Denote by $\Pi_n$
 the set of all compositions of $n$ which are peak compositions.
  {}From the above discussion on sets of peaks, it is obvious that
\begin{equation}
\label{eq:DescripPeakComp}%
\Pi_n=\set{\mathbf c=(c_1,\dots,c_k)\vDash n \st
c_{i}>1,\  2\le i\le k-1}.
\end{equation}
Moreover, it is clear
 that for all $\sigma,\tau\in\SG_n$ we have
$\ipeak(\sigma)\subseteq\peak(\sigma)\subseteq\Des(\sigma)$
 and  $D(\sigma)=D(\tau)\implies P(\sigma)=
P(\tau)\implies \IP(\sigma)=\IP(\tau)$.

 For each $\mathbf c\in\Pi_n$, we set
$
p_{\mathbf c} =\sum_{P(\sigma)=\mathbf c} \sigma.
$
Then Aguiar, Bergeron and Nyman \cite{aguiar bergeron nyman}
 have shown that
$$
\Peak_n=\bigoplus_{\mathbf c\in\Pi_n} \ZM p_{\mathbf c}
$$
is a subalgebra of $\Sigma_n$. We call $\Peak_n$ the {\sl exterior
peak algebra}. They have also shown that the graded space
$\Peak=\mathcal E(P)$ is a subcoalgebra of $(\k[\SG],\Delta)$.

In \cite{schocker}, the author has given an analog of Atkinson's
proof of Solomon's result for the  peak algebra, see
\S\ref{ss:IntPeaks}. Here we give a similar proof for the exterior
peak algebra.  For $s_i\in S_n$ and $\sigma\in\SG_n$, we
define the exterior peak equivalence relation as follows:
\begin{eqnarray}
\label{eq:EquivPeakSn} \sigma \dot\sim_P s_i\sigma &\iff&
\left\{\begin{array}{l}%
s_i\sigma\dot\sim_D \sigma\quad
\textrm{or}\\
i=1\quad \textrm{ and }\quad s_1\sigma\not = \sigma s_1,
\end{array}\right.\\
&\iff& \left\{\begin{array}{l}%
|\sigma^{-1}(i)-\sigma^{-1}(i+1)|>1
\quad \textrm{or}\\
i=1\quad \textrm{and }\quad  [\sigma(1),\sigma(2)]\not = [1,2].
\end{array}\right.
\end{eqnarray}
It is easily seen that $\dot\sim_P$ is a symmetric relation. The
reflexive and transitive closure of $\dot\sim_P$ is called
\textit{the exterior peak equivalence}, and is denoted by
$\sim_P$.

\begin{lem}
\label{le:PeakSets} For $\sigma\in\SG_n$ and $s_i\in S_n$, if
$s_i\sigma\dot\sim_P \sigma$, then $P(s_i\sigma)=P(\sigma)$.
\end{lem}
\begin{proof} By
Proposition~\ref{pr:AtkinsonDescent} and Eq.~(\ref{eq:EquivPeakSn}) we
can easily reduce the proof to the case where $i=1$, the letters
$1,2$ are adjacent and $[\sigma(1),\sigma(2)]\not = [1,2]$. As the
word $1\cdot2$ (or $2\cdot 1$) is not a prefix of the word
$\sigma$, it is obvious that exchanging $1$ and $2$ does not
change $\peak(\sigma)$.
\end{proof}

\begin{prop}
\label{pr:PeakSetsEquiva} For $\sigma,\tau\in\SG_n$,
$P(\sigma)=P(\tau)$ if and only if $\sigma \sim_P \tau$.
\end{prop}
\begin{proof} By Lemma~\ref{le:PeakSets} we have only to show that
if $P(\sigma)=P(\tau)$ then  $\sigma \sim_P \tau$.

Assume that $I=\Des(\sigma)\not =\peak(\sigma)$, then there is $j\in I$
such that $j,j-1\in I$ and $j+1\notin I$. Let
$\sigma'\in\SG_n$ be such that $\Des(\sigma')=I$ and such that
$\sigma'(j)=2=\sigma'(j+1)+1$ (Corollary~\ref{co:DescentSets}).
We have $s_1\sigma'\dot\sim_P \sigma'$, that is
$P(\sigma)=P(s_1\sigma')$ and
$|\Des(s_1\sigma')|=|\Des(\sigma')|-1$. By induction on $\Des(\sigma)$
and by Proposition~\ref{pr:AtkinsonDescent} (ii) we obtain an
element $\sigma_1\in \SG_n$ such that $\sigma\sim_D \sigma'\sim_P
\sigma_1$ and $\Des(\sigma_1)=\peak(\sigma_1)=\peak(\sigma)$.
In other words $\sigma\sim_P \sigma_1$
since  two permutations having the same descent sets have the
 same exterior peak set. Proceed similarly with $\tau$  to obtain a permutation
 $\tau_1$ such that $\Des(\tau_1)=\peak(\tau_1)=\peak(\tau)$. Then $\tau\sim_P\tau_1$.
 As $\Des(\tau_1)=\Des(\sigma_1)$ and as $\tau_1\sim_D\sigma_1$ implies
 $\tau_1\sim_P\sigma_1$ the proposition follows from Proposition~\ref{pr:AtkinsonDescent}.
\end{proof}

The following corollary is immediate.
\begin{cor}
\label{co:PeakLeftConnected} Each equivalence class under $\sim_P$
is left-connected.
\end{cor}

If $\sigma\in\SG_n$ and $\mathbf c,\mathbf d\in\Pi_n$, we set
$$
P(\mathbf c,\mathbf d,\sigma) = \set{(u,v) \in \SG_n \times \SG_n
\st P(u)=\mathbf c,\ P(v)=\mathbf d,\ uv=\sigma }.
$$
For $s_i \in S_n$ and  $(u,v) \in \SG_n \times \SG_n$, we define
$$
\psi_i^P (u,v) = \left\{\begin{array}{ll}
\psi^D_i (u,v)&\textrm{if } s_i u\dot\sim_D u\\
(s_1 u,v)&\textrm{if } i=1\ \textrm{and}\ s_1u\not= us_1\\
(u,s_1 v)&\textrm{if } i=1\ \textrm{and}\ s_1u = us_1\,.
\end{array}\right.
$$
We have $(\psi^P_i)^2 = \text{Id}_{\SG_n \times \SG_n}$. In particular
$\psi^P_i$ is a bijection.

\begin{lem}
\label{le:PsiPeak} If $\sigma \in \SG_n$ and $s_i \in S_n$ are
such that $\sigma \dot\sim_P s_i \sigma$, then
$$\psi^P_i\big(P(\mathbf c,\mathbf d,\sigma)\big)=P(\mathbf
c,\mathbf d,s_i\sigma),\quad\text{for all }\mathbf c,\mathbf d \in
\Pi_n.$$
\end{lem}
\begin{proof} Let $(u,v)\in P(\mathbf c,\mathbf d,\sigma)$.
 First assume that $\sigma \dot\sim_P s_1 \sigma$. Observe that if $s_1u=us_1$ and $s_1v=vs_1$,
 then $s_1\sigma = \sigma s_1$ which is a contradiction. In this case  $\psi_i^P
 (u,v)\in P(\mathbf c,\mathbf d,s_i\sigma)$.
 Assume next that $\sigma \dot\sim_D s_i \sigma$. We observe from Lemma~\ref{le:Psi} and $\peak(\sigma)\subseteq \Des(\sigma)$ that if $s_i
u\dot\sim_D u$, then $\psi_i^P (u,v)=\psi_i^D (u,v)\in
\psi^D_i\big(D(D(u),D(v),\sigma)\big)=D\big(D(u),D(v),s_i\sigma\big)\subseteq P(\mathbf c,\mathbf d,s_i\sigma)$. \end{proof}

\begin{lem}
\label{le:PeakCoalgebra}%
The graded connected map $P$ has the restriction property.
\end{lem}
\begin{proof} Proceed as in the proof of Lemma~\ref{le:EquivHaveRP}.
\end{proof}

\begin{rem}\label{peakInduction} \textnormal{The graded map $P$ does not have the induction
property. For example, the set
$X_{(1,2)}(\set{1}\times\set{12})=\set{123,213,312}$ does not
contain $321\dot\sim_P 312$.}
\end{rem}

Then Lemmas~\ref{le:PsiPeak} and \ref{le:PeakCoalgebra},
Proposition~\ref{pr:PeakSetsEquiva} and Theorem~\ref{th:RP} imply
directly the following theorem.

\begin{thm}
\label{th:PeakAlgebra} %
$\quad$
\begin{enumerate}

\item[(i)] $\Peak_n$ is a subalgebra of $\Sigma_n$.

\item[(ii)] $\Peak$ is a subcoalgebra of $(\k[\SG],\Delta)$.

\end{enumerate}
\end{thm}

In general the graded connected map $P^G$ defined by
Eq.~(\ref{eq:DefRhoG}) does not have the induction property and nor the
restriction property.

\subsection{Coalgebras in the hyperoctahedral group}
\label{ss:HyperGroup}

We give here some examples of coalgebras associated to
 graded map which are not induced from graded maps on the symmetric group.

For $G=\ZM/2\ZM=\set{-1,+1}$,
the hyperoctahedral group $G_n$ is a Coxeter group of type
$B_n$ generated
 by $S_n^\pm=\set{s_0,s_1,\dots,s_{n-1}}=\set{s_0}\cup S_n$ where $s_0= 1_n^{(-1,+1,+1,\dots,+1)}$.
The Coxeter length of $\alpha\in G_n$ is in this case
$$
\ell_B(\alpha)=\min\set{k>0\st \alpha=r_1\dots r_k,\ r_i\in S^\pm_n}.
$$
  Recall that a subset $X$ of
$G_n$ is \textit{left}-B-\textit{connected} if for each $\alpha,\beta\in X$,
there is a sequence $\alpha = \alpha_1 , \alpha_2,\dots ,\alpha_k
= \beta$ of elements in $X$ such that
$\alpha_{j+1}\alpha_{j}^{-1}\in S^\pm_n$, for all $j\in[1,k-1]$. That
is, these classes can be seen as a set of adjacent nodes in the
{\sl type} B-\textit{permutahedron}.

 The B-\textit{descent set} of $\alpha\in G_n$, with the convention that $\alpha(0)=0$, is
\begin{eqnarray}
\label{eq:DefDescB}
\Des_B(\alpha)&=&\set{i\in[0,n-1]\st\ell_B(\alpha s_i)<\ell_B(\alpha)} \\
&=&\set{i\in[0,n-1]\st \alpha(i)>\alpha(i+1)}.\nonumber
\end{eqnarray}
We get a graded map $\Des_B\colon\mathsf G\to\bigoplus_{n\in \NM}[0,n-1]$ and $\mathcal E(\Des_B)_n$ is precisely the {\sl Solomon
descent algebra} $\Sigma(G_n)$ associated to $G_n$ \cite{solomon}.
Atkinson \cite{atkinson} has shown an analog of
Proposition~\ref{pr:AtkinsonDescent}
 for B-descent sets: B-descent sets are left-B-connected. More
precisely, the type B descent equivalence is defined as
 the transitive and reflexive closure of the following relation:
 if $r\in S^\pm_n$ and $\alpha\in G_n$, then
\begin{equation*}
\label{eq:DefDescEquivB}
r\alpha\dot\sim_{\Des_B} \alpha \iff \alpha^{-1}
r\alpha\not\in S^\pm_n ,
\end{equation*}
or equivalently
\begin{equation*}
\label{eq:DefDescEquivBBIS}
r\alpha\dot\sim_{\Des_B} \alpha \iff
\left\{\begin{array}{l}
r\in S_n \quad\textrm{and}\quad r\alpha\dot\sim_D^G \alpha\\
r=s_0\quad\textrm{and}\quad |\alpha(1)|>1
\end{array}\right.
\end{equation*}
since  $s_0$ is not in the conjugacy class of $s_i\in S_n$.
 Obviously  $\alpha\sim_{D}^G \beta\implies\alpha\sim_{\Des_B}\beta$. Hence
 $\Sigma(G_n)\subseteq \Sigma_n(G)$.

In \cite{aguiar bergeron nyman,bmws}, the authors have shown that
$\Sigma_B=\mathcal E(\Des_B)$ is a subcoalgebra of $(\Sigma(G),\Delta)$.
 With our theory we obtain this result as a consequence of
 the following lemma.

\begin{lem}
\label{le:TypeBDescent}
The graded map $\Des_B$ has the restriction property.
\end{lem}
\begin{proof} Let $n,m\in\NM$ and $u\in X_{(n,m)}$.
 For $\alpha\in G_n$ and $\beta\in G_m$,
we have to show that for $r\in S_n^\pm$, either of the following holds:

\item{(i)} if $r\alpha\dot\sim_{\Des_B} \alpha$,
 then $(r\alpha\times \beta)u^{-1}\dot\sim_{\Des_B}
(\alpha\times \beta)u^{-1}$,

\item{(ii)}  if $r\beta\dot\sim_{\Des_B} \beta$,
 then $\Des_B\big((\alpha\times r\beta)u^{-1})=\Des_B\big(
(\alpha\times \beta)u^{-1}\big)$.

\noindent
 By  Eq.~(\ref{eq:DefDescEquivBBIS}) we have to consider only $r=s_0$,
 since  Theorem~\ref{th:Generating} gives us that  $D_G$ has the restriction property.
It is well known that $\SG_{(n.m)}$ is a parabolic subgroup of $G_{n+m}$ and that
 $\ell_B((v_1\times v_2) u^{-1})=\ell_B(v_1\times v_2)+\ell_B( u^{-1})$.
The lemma follows easily using Eq.~(\ref{eq:DefDescB}).
\end{proof}

\begin{rem} \textnormal{
In symmetric groups  the graded map $\mathcal T$ gives us a construction
of the  {\sl associahedron} from the permutohedron (see for
instance \cite{loday ronco}). For the hyperoctahedral groups,
 such a map has been  described by
 Reiner \cite{reiner}. Let us denote by $\mathsf{Gon}_n$ the set of {\sl centrally
symmetric $(2n+2)$-gons} and let
$$
\mathcal G\colon\mathsf G\to  \mathsf{Gon}=\bigoplus_{n\in\NM}\mathsf{Gon}_n
$$
be the graded map given by Reiner. This map  gives a construction
of the {\sl cyclohedron} from the type $B$-permutohedron. It is a
left-B-connected graded map.
 In a work in progress
\cite{aguiar thomas}, Aguiar and Thomas have observed that $\mathcal E(\mathcal G)$ is
 a subcoalgebra of $(\k[\mathsf G],\Delta)$.
Recently, Reading \cite{reading} has defined a $G$-equivalence relation
 whose  equivalence classes are parameterized by the centrally symmetric
$(2n+2)$-gons. We could prove using this relation that $\mathcal G$ has
 the restriction property, and then, we could obtain another proof
 of the above result of Aguiar and Thomas. It is interesting to notice
 that the number of  centrally symmetric $(2n+2)$-gons is
$\left(\begin{array}{@{}c@{}}2n\\n\end{array}\right)$ which is also the number of $G$-sequences
 of trees whose sum of vertices is $n$ (see Eq.~(\ref{eq:CardinalInducedEquiv})).
 That means that the dimension of
 homogeneous parts of $\mathcal E(\mathcal G)$ and $\mathcal E(\mathcal T_G)$
 are equal.}
\end{rem}

%
%

\section{The $G$-descents to $G$-peaks map}
\label{se:Duality}

\subsection{The $\Theta_G$-function}
\label{ss:Theta}

There is a well known Hopf endomorphism $\Theta$ on symmetric
functions whose image is the space spanned by $Q$-Schur functions
\cite{macdonald}. Several authors \cite{aguiar bergeron nyman,
aguiar bergeron sottile, bht, stembridge} have extended this
morphisms to the quasi-symmetric functions, Hopf algebras and
$(\Sigma,*,\Delta)$. It plays an important role in studying the
peak algebra. We first recall the definition of
$\Theta\colon\Sigma\to \Sigma$.

We have that $\Sigma= {\mathcal E}(D)=\bigoplus_{n\ge 0}
\k\otimes_{\ZM}\Sigma_n$, where $\Sigma_n$ is the free
$\ZM$-module spanned by $\set{d_{\mathbf c}\st {\mathbf c}\models
n}$  (see \S\ref{ss:DescentSets}). It is well known that $\Sigma$
is freely generated as an algebra by $\set{d_{(n)}=1_n\st n\ge 1}$
(see \cite{malvenuto reutenauer, loday ronco, thibon}). We set
\begin{equation}\label{eq:defTheta}
 \Theta(d_{(n)})=2 \ip_{(n)}.
\end{equation}
This defines a unique morphism of algebras, and it is
straightforward to check that
$\Theta(\Delta(d_{(n)}))=\Delta(\ip_{(n)})$. Hence $\Theta\colon
\Sigma\to\Sigma$ is a Hopf morphism whose image is $\IPeak$. It is
exactly the morphism defined in \cite{aguiar bergeron nyman,
aguiar bergeron sottile, bht, stembridge}.

For any ${\mathbf c}\models n$, there is an explicit formula for
$\Theta(d_{\mathbf c})$ (see for example \cite[Thm~5.8]{aguiar
bergeron nyman}. For this we need to introduce more notation on
compositions. For ${\mathbf c}=(c_1,c_2,\ldots c_s)\models n$ we
denote the {\sl number of parts} of $\mathbf c$ by $\tau({\mathbf
c})=s$.  There is a natural map from compositions of $n$ to peak
compositions. For ${\mathbf c}=(c_1,c_2,\ldots c_s)\models n$, let
$\Lambda({\mathbf c})\in\IPi_n$ be the peak composition obtained
from $\mathbf c$ as follows. The composition $\mathbf c$
factorizes uniquely into compositions ${\mathbf c}_1\cdot {\mathbf
c}_2\cdots {\mathbf c}_r$ where for $i<r$ we have ${\mathbf
c}_i=(1,1,\ldots,1,m)$ and $m>1$, and ${\mathbf
c}_r=(1,1,\ldots,1)$. In this factorization, any sequence of $1$'s
may be empty. We then define
$$
\Lambda({\mathbf c})=(\|{\mathbf c}_1\|,\|{\mathbf
c}_2\|,\ldots,\|{\mathbf c}_r\|)\in\IPi_n,
$$
removing the last part if it is zero. For example,
 $$\Lambda\big((1,2,1,1,3,3)\big)=\big(\|(1,2)\|,\|(1,1,3)\|,\|(3)\|\big)=(3,5,3).$$
Similarly,  $\Lambda\big((2,1,3,1)\big)=(2,4,1)$, and
$\Lambda\big((2,1,1)\big)=(2,2)$. Finally, recall from
Eq.~(\ref{eq:BijCompSubset1}) that $I_{\mathbf c}$ is a subset of
$[1,n-1]$. For any ${\mathbf c}\models n$,
\begin{equation}\label{eq:Thetac}
 \Theta(d_{\mathbf c})=\sum_{{\mathbf e}\in\Phi(\mathbf c)} 2^{\tau(\Lambda({\mathbf e}))} d_{\mathbf e},
\end{equation}
where
\begin{equation}\label{eq:sumPhi}
\Phi({\mathbf c})=\set{{\mathbf e}\models \|\mathbf c\|\, \big|\,\big( i\in I_{\Lambda({\mathbf e})}\  \implies\   \big| \set{i-1,i}\cap I_{\mathbf c}\big|=1\big)}.
\end{equation}

Since $D$ satisfies (FP) and (IP), Corollary~\ref{co:FreeIP}
and~\ref{co:Fonctoriality} give that the assignment
$G\rightsquigarrow  {\mathcal E}(D_G)$ is functorial and
$\Sigma(G)= {\mathcal E}(D_G)$ is freely generated as an algebra
by
  $$\set{d^G_{(n)^g}=(1_n)^g\st n\ge 1,\ g\in G}.$$
\begin{prop}
\label{pr:Theta} The map $\Theta_G\colon\Sigma(G)\to\Sigma(G)$ defined by
\begin{equation}\label{eq:defThetaG}
 \Theta(d^G_{(n)^g})=\ip^G_{(n)^g}
\end{equation}
is a Hopf morphism. Moreover,
$\Theta_G\big(\Sigma(G)\big)=\IPeak(G)$.
\end{prop}
\begin{proof}
Eq.~(\ref{eq:defThetaG}) defines a unique morphism of algebras. It
is straightforward to check that
$\Theta(\Delta(d^G_{(n)^g}))=\Delta(\ip^G_{(n)^g})$. Hence
$\Theta_G$ is a Hopf morphism whose image, by
Theorem~\ref{th:PeakHopf}, is $\IP(G)$.
\end{proof}

For any $G$-composition ${\mathbf c}^{\mathbf g}$, we are
interested in an explicit formula for $\Theta_G(d_{{\mathbf
c}^{\mathbf g}}^G)$. For this, let ${\mathbf c}^{\mathbf
g}={\mathbf c_{1}}^{g_1}\cdot {\mathbf c_{1}}^{ g_2} \cdots
{\mathbf c_{k}}^{ g_k}$ be the rainbow decomposition of ${\mathbf
c}^{\mathbf g}$. {}From Corollary~\ref{co:FreeIP} we clearly have
that
$$
d_{{\mathbf c}^{\mathbf g}}^G = d_{{\mathbf c_{1}}^{g_1}}^G
*\cdots * d_{{\mathbf c_{k}}^{g_k}}^G
          = \big(d_{{\mathbf c_{1}}}^G\big)^{g_1} *\cdots * \big(d_{{\mathbf c_{k}}}^G\big)^{g_k}.
$$
Now, from Eq.~(\ref{eq:defThetaG}), it is also clear that for a single $g\in G$, we have
 $$\Theta(d_{{\mathbf c}^{g}}^G)=\Theta\big((d_{{\mathbf c}}^G)^{g}\big)=\big(\Theta(d_{{\mathbf c}}^G)\big)^{g}.$$
 Combining this with Eq.~(\ref{eq:Thetac}) we obtain the following theorem.

\begin{thm}
Let $G$ be a finite abelian group and let ${\mathbf c}^{\mathbf
g}={\mathbf c_{1}}^{g_1}\cdot {\mathbf c_{2}}^{ g_2} \cdots
{\mathbf c_{k}}^{ g_k}$ be the rainbow decomposition of the
$G$-composition ${\mathbf c}^{\mathbf g}$. We have
\begin{equation}\label{eq:Thetacg}
 \Theta_G (d^G_{{\mathbf c}^\mathbf{g}})=
    \sum_{1\le i\le k,\, {\mathbf e}_{i}\in\Phi(\mathbf c_{i})}
    2^{\sum_{i=1}^k \tau(\Lambda({\mathbf e}_{i}))} \,
    d^G_{{\mathbf e_{1}}^{g_1}\cdot {\mathbf e_{2}}^{g_2} \cdots {\mathbf e_{k}}^{ g_k}}\,.
\end{equation}
\end{thm}

\subsection{The dual side}
\label{ss:DualSide}

Let $\mathcal A$ be a countable alphabet viewed as a set
of variables and $G$ be a finite abelian group. We let $\mathcal Z=\mathcal A^G$ be the alphabet where we put an arbitrary total order on $G$ and we order $\mathcal Z$ as follows. For $a^g,b^h\in\mathcal Z$, we set $a^g<b^h$ if and only if ($a<b$) or ($a=b$ and $g<h$).
In this section we consider subspaces of  the space $\k[[\mathcal Z]]$ of formal series in the variables
$a^g\in\mathcal Z$.

 Let $\mathbf c^\mathbf{g}=(c_1^{g_1},c_2^{g_2},\ldots,c_k^{g_k})$ be a
$G$-composition of  $n$ and  set $t_i=c_1+c_2+\cdots+c_i$ for each
$i$. We denote by
$$(h_1,h_2,\ldots,h_n)=(\underbrace{g_1,g_1,\ldots,
g_1}_{\text{$c_1$
times}},\underbrace{g_2,g_2,\ldots,g_2}_{\text{$c_2$
times}},\ldots,\underbrace{g_k,g_k,\ldots,g_k}_{\text{$c_k$
times}}).$$ We let
$\mathcal Z_{\mathbf c^\mathbf{g}}$ be the set of all $n$-tuples
$(z_1,z_2,\ldots,z_n)\in \mathcal Z^n$ where $z_i=a_i^{h_i}$ for some $a_1\le a_2\le\cdots\le a_n$ such that
$$\forall i\in\{1,2,\ldots,k-1\},\quad g_i\geq g_{i+1}\Longrightarrow
a_{t_i}<a_{t_i+1}.$$
Finally, we define the formal series in $\k[[\mathcal Z]]$
$$
F_{\mathbf c^\mathbf{g}}=\sum_{(z_1,z_2,\ldots,z_n)\in {\mathcal Z}_{\mathbf
c^\mathbf{g}}}z_1z_2\cdots z_n.
$$

In \cite{baumann hohlweg} the authors have shown that the vector
subspace $\mathrm{QSym}(G)$ spanned by the
$F_\mathbf{c^\mathbf{g}}$ for all $G$-compositions is a subalgebra
of $\k[[\mathcal A^G]]$. It is a graded connected Hopf algebra
that is graded dual to $\Sigma(G)$. In particular, the basis
$\set{F_{\mathbf c^\mathbf{g}}}$ of $\QSym(G)$ is dual to the
basis $\set{d^G_{{\mathbf c}^\mathbf{g}}}$ of $\Sigma(G)$. To see
this, we first recall that $\k[G]$ is self dual. If $\alpha\in
G_n$, then $\alpha^{-1}=\alpha^*$ is its dual element in the dual
basis \cite[\S5.1]{baumann hohlweg}. For any $G$-composition
$\mathbf c$, \  $\alpha^* (d^G_{\mathbf c}) = 1$ iff
 $D_G(\alpha^{-1})=\mathbf c$.
But by \cite[Prop.~22 and Diag.~(11)]{baumann hohlweg} we know that
 the dual of $\Sigma(G$) is obtained as the quotient  of $\k[G]$ by the
relation $\alpha \sim_{R_G} \beta$ where
$R_G(\alpha)=D_G(\alpha^{-1})$. For $\alpha\in G_n$ such that
$R_G(\alpha)=\mathbf c$ we have from \cite[Thm~33 and
Prop.~34]{baumann hohlweg} that the image of $\alpha^*$ in this
quotient correspond to $F_{\mathbf c}$ of $\QSym(G)$. We then have
that $F_{\mathbf c}(d_{\mathbf e}^G)= 1$ if and only if $\mathbf
c=\mathbf e$.

For a $G$-composition ${\mathbf c}^{\mathbf g}$ let ${\mathbf
c}^{\mathbf g}={\mathbf c_{(1)}}^{g_1}\cdot {\mathbf c_{(2)}}^{
g_2} \cdots {\mathbf c_{(k)}}^{ g_k}$ be its rainbow decomposition
and let $n_i=\|{\mathbf c}_{(i)}\|$. Using Eq.~(\ref{eq:Thetacg}),
the morphism $\Theta_G$ has a graded dual $\Theta_G^*\colon
\QSym(G)\to \QSym(G)$ given by
\begin{equation}\label{eq:ThetaDual}
 \Theta^*_G (F_{{\mathbf c}^\mathbf{g}})=2^{\sum_{i=1}^k \tau(\Lambda({\mathbf c}_{i}))}
    \sum_{1\le i\le k,\,{\mathbf e}_{i}\in \Phi^*({\mathbf c}_{i})}
    F_{{\mathbf e_{1}}^{g_1}\cdot {\mathbf e_{2}}^{g_2} \cdots {\mathbf e_{k}}^{ g_k}},
\end{equation}
where
\begin{equation}\label{eq:sumPhidual}
\Phi^*({\mathbf c})=\set{{\mathbf e}\models \|\mathbf c\|\, \big|\,\big( i\in I_{\Lambda({\mathbf c})}\  \implies\   \big| \set{i-1,i}\cap I_{\mathbf e}\big|=1\big)}.
\end{equation}

We deduce the following theorem

\begin{thm}\label{th:DualPeak} The image of $\Theta^*_G$ is the
graded dual Hopf algebra $\IPeak(G)^*$ of $\IPeak(G)$.
\end{thm}

\subsection{$G$-coloured combinatorial Hopf algebras}
\label{ss:CombinHopfAlg}

In \cite{aguiar bergeron sottile} we introduced the theory of
combinatorial Hopf algebras and associated objects. The framework
consists of pairs $({\mathcal H},\zeta)$ where $\mathcal H$ is a
connected graded Hopf algebra and $\zeta\colon\mathcal H\to\k$ is
an algebra morphism (called there {\sl character}). We have shown
in that paper that for a certain character $\zeta_{\mathcal Q}$
the pair $(\QSym,\zeta_{\mathcal Q})$ is the terminal object for
the category of pairs $({\mathcal H},\zeta)$. Moreover we have
shown that $\IPeak$ was the so-called odd subalgebra of
$(\QSym,\zeta_{\mathcal Q})$. These algebras play an important
role among the connected graded Hopf algebras and we refer the
reader to \cite{aguiar bergeron sottile} for more motivation. Here
we are interested in a similar theory for the $G$-coloured
version. We will just outline the ideas, as this should be part of
forthcoming work.

For $G=\set{1^G}$, The map $\zeta_{\mathcal Q}\colon \QSym\to\k$
is the morphism of algebras defined by $\zeta_{\mathcal
Q}(F)=F(1,0,0,\cdots)$. The functoriality of the construction of
$\QSym(G)$ implies that the following diagram commutes

\vskip-50pt
\begin{equation}\label{di:Theta}
    \raise -39pt\hbox{
       \begin{picture}(140,158)
       \put(-5,90){$\mathrm{QSym}(G)$}    \put(85,90){$\mathrm{QSym}(G)$}
       \put(105,70){$\varphi$} \put(10,70){$\varphi$}
       \put(0,40){$\mathrm{QSym}$}    \put(88,40){$\mathrm{QSym}$}
       \put(58,0){$\k$}
       \put(53,100){$\Theta_G^*$}
       \put(53,50){$\Theta_\set{1}^*$}
       \put(24,16){$\nu$}
       \put(88,16){$\zeta_{\mathcal Q}$}
       \put(20,85){\vector(0,-1){35}}
       \put(100,85){\vector(0,-1){35}}
       \put(40,43){\vector(1,0){40}}
       \put(40,93){\vector(1,0){40}}
       \put(20,35){\vector(1,-1){27}}
       \put(100,35){\vector(-1,-1){27}}
      \end{picture}} ,
\end{equation}
where $\nu=\zeta_{\mathcal Q}\circ \Theta^*_\set{1^G}$ and
$\varphi$ is the Hopf morphism induced by the inclusion
$\set{1^G}\to G$. We have shown in \cite[Prop.~6.4]{aguiar
bergeron sottile} that  the pair $(QSym(\set{1}),\nu)$ is odd,
hence $(\QSym(G),\nu\circ\varphi)$ is odd. Using
\cite[Prop.~6.1]{aguiar bergeron sottile} in such situation gives
us that the image of $\Theta_G^*$, that is $\IPeak(G)^*$, must be
contained in the odd subalgebra of $\QSym(G)$.  But a simple
dimension count gives us that the odd subalgebra of
$(\QSym(G),\zeta_{\mathcal Q}\circ\varphi)$ is larger than
$\IPeak(G)^*$. On the other hand, if we restrict ourselves on a
smaller category, then an interesting theory unfolds.

Let $G$ be a fixed finite abelian group. Consider pairs
$({\mathcal H},\zeta^G)$ where $\mathcal H$ is a connected graded
Hopf algebra and $\zeta\colon\mathcal H\to\k[G]$ is an algebra
morphism. In analogy with \cite{aguiar bergeron sottile} we will
call such a morphism $\zeta\colon\mathcal H\to\k[G]$ a $G$-{\sl
character}. The category of combinatorial $G$-coloured Hopf
algebras consists of pairs $({\mathcal H},\zeta^G)$ as above and
graded Hopf morphisms $\Psi$ such that
 \vskip-50pt
$$    \raise -39pt\hbox{
       \begin{picture}(140,108)
       \put(25,40){$\mathcal H$}    \put(86,40){$\mathcal H'$}
       \put(52,0){$\k[G]$}
       \put(53,50){$\Psi$}
       \put(28,18){$\zeta^G$}
       \put(86,18){$\zeta'^G$}
       \put(40,43){\vector(1,0){40}}
       \put(32,35){\vector(1,-1){22}}
       \put(86,35){\vector(-1,-1){22}}
      \end{picture}}
$$
commutes.

There is a natural algebra morphism $\k[\mathcal A^G]\to\k[G]$
which maps the variable $a_i^g\in \mathcal A^G$ to $g$ if $i=1$
and to zero otherwise. This gives us the algebra morphism
$\zeta_{\mathcal Q}^G\colon \QSym(G)\to \QSym(G)$ where
  \begin{equation}\label{eq:zeta0}
    \zeta_{\mathcal Q}^G(F_{\mathbf c^{\mathbf g}})\ =\ \left\{
                   \begin{array}{ll}g&\textrm{if}\
              \mathbf c=(n)\textrm{ and }\mathbf g=g,\\
               1&\textrm{if}\  \mathbf c=(),\\
               0&\textrm{otherwise}.\end{array}
   \right.  \end{equation}
We have a theorem very similar to \cite[Thm.~4.1]{aguiar bergeron
sottile} with essentially  the same proof.

\begin{thm}
The pair $(\QSym(G),\zeta_{\mathcal Q}^G)$ is the terminal object
in the category of combinatorial $G$-coloured Hopf algebras.
\end{thm}

We then define the convolution product of two $G$-characters as
follow. For two $G$-characters $\zeta^G,\nu^G\colon\mathcal
H\to\k[G]$ let
 $$ \zeta^G * \nu^G= m_{\k[G]}\circ (\zeta^G\otimes \nu^G)\circ\Delta_{\mathcal H}.
 $$
 $G$-characters form a group with inverse (for the convolution) given by
 $(\zeta^G)^{-1}=\zeta^G\circ S_{\mathcal H}$.  For a $G$-character
 $\zeta^G \colon\mathcal H\to\k[G]$ we let $\overline{\zeta^G} \colon\mathcal H\to\k[G]$
be the character such that for a homogeneous element $h\in\mathcal H$ of
degree $n$,\quad $\overline{\zeta^G}(h)=(-1)^n\zeta^G(h)$. The odd $G$-subalgebra
$S_-({\mathcal H},\zeta^G)$ of a pair $({\mathcal H},\zeta^G)$ is defined as in
\cite[Def.~5.7]{aguiar bergeron sottile}, namely the largest graded subcoalgebra of
 $\mathcal H$ such that for all $h\in S_-({\mathcal H},\zeta^G)$ we have $\overline{\zeta^G}(h)=(\zeta^G)^{-1}(h)$.
It is in this context that
 \begin{equation}\label{eq:TheEnd}
 S_-(\QSym(G),\zeta^G_{\mathcal Q}) = \IP(G)^*.
\end{equation}
Combining this with an analog of \cite[Cor.~6.2]{aguiar bergeron
sottile}, we have that $(\IPeak(G)^*,\zeta_{\mathcal Q}^G)$ is the
terminal object of odd combinatorial $G$-coloured Hopf algebras.

%
%

\end{document}